\documentclass{article}

\usepackage[english]{babel}
\usepackage{amsmath}
\usepackage{graphicx}
\usepackage[colorlinks=true, allcolors=blue]{hyperref}
\usepackage{bbm}
\usepackage{amssymb}
\usepackage{xcolor}
\usepackage{amsthm}
\usepackage{caption}
\usepackage{subcaption}
\usepackage{cleveref}
\usepackage{booktabs}
\usepackage{tabularx}
\usepackage{authblk}
\usepackage{doi}
\usepackage[normalem]{ulem}
\usepackage{setspace}
\usepackage[authoryear]{natbib}
\usepackage{enumitem}
\setlist[itemize]{topsep=1pt, itemsep=1pt, parsep=1pt}
\usepackage{tikz}

\usepackage{placeins}

\usepackage[letterpaper,top=2cm,bottom=2cm,left=3cm,right=3cm,marginparwidth=1.75cm]{geometry}

\theoremstyle{definition}

\newtheorem{ex}{Example}

\newcommand{\new}[1]{\textcolor{black}{#1}}

\renewcommand{\sout}[1]{}

\title{\new{Non-monotone clustering and energetic regime switching in aggregation–diffusion systems}}

\author[1,2]{Nathalie Wehlitz}
\author[3]{Georg Richard Scherzer}
\author[4]{Carsten Hartmann}
\author[1]{Stefanie Winkelmann\footnote{Email of corresponding author: \href{mailto:winkelmann@zib.de}{winkelmann@zib.de}}}

\affil[1]{Zuse Institute Berlin, Berlin, Germany}
\affil[2]{Freie Universität Berlin, Berlin, Germany}
\affil[3]{Technische Universität Berlin, Berlin, Germany}
\affil[4]{Brandenburgische Technische Universität, Cottbus--Senftenberg, Germany}

\date{}

\begin{document}

\maketitle

\doublespacing

\begin{abstract}
We investigate transient clustering dynamics in nonlocal aggregation--diffusion systems from an energetic perspective. Starting from a stochastic interacting particle system, we study the associated macroscopic McKean--Vlasov equation on the torus and \new{use} its Wasserstein gradient-flow structure to analyse the competition between interaction-driven aggregation and entropy-driven diffusion. Through numerical experiments for locally attractive interaction kernels, \new{we observe non-monotone clustering behaviour in intuitive geometric
observables: in particular, the density peak height may first increase and
subsequently decrease along trajectories converging to equilibrium.
Using the decomposition of the free energy as a diagnostic, we identify
alternating aggregation- and diffusion-dominated regimes underlying these
transient dynamics. 
We investigate to what extent transitions between these regimes are related to the evolution of geometric observables, such as the peak height. Our results illustrate how non-monotone
transient clustering can arise even in a dissipative gradient-flow system and
highlight the complementary information provided by 
geometric clustering observables and the free-energy-based regime classification.}

\quad

\end{abstract}

\section{Introduction}

Collective dynamics and clustering phenomena arise in a wide range of interacting particle systems, including flocking or swarming dynamics~\citep{carrillo2010, topaz2008} and opinion formation processes~\citep{garnier2017consensus, helfmann2021, zimper2026}, where clustering corresponds to the emergence of consensus. Such systems are often modelled through locally attractive interactions combined with stochastic effects, leading to rich pattern formation and complex relaxation behaviour.

The long-time behaviour of aggregation--diffusion systems has been extensively studied in recent years. Current research focuses include the characterisation of stationary states, bifurcation structure, phase transitions, coarse-graining limits, and variational formulations of the dynamics~\citep{carrillo2019, carrillo_long-time_2020, adams2026}, including analyses of equilibrium states and their stability properties for noisy nonlocal aggregation models~\citep{yang2026noisy}. In particular, for suitable interaction potentials and sufficiently small noise strengths, the associated free energy becomes nonconvex, leading to the coexistence of multiple locally stable equilibria, including both homogeneous and clustered stationary states~\citep{carrillo_long-time_2020,chayes2010mckean}. In this coexistence regime, the asymptotic state depends sensitively on the initial condition. Recent work has further connected clustering behaviour with discontinuous phase transitions and variational properties of the free-energy landscape~\citep{gerber2025formation}. Related clustering phenomena and metastable behaviour have also been investigated for kinetic and underdamped McKean--Vlasov dynamics~\citep{leimkuhler2026cluster}.

Most existing works focus primarily on equilibrium configurations, long-time convergence, or linear stability properties. In contrast, transient dynamics remain much less understood. In particular, the coexistence of multiple equilibria is a static property and does not explain the transient clustering behaviour observed along trajectories approaching equilibrium. The present work focuses on these transient dynamics and complements recent work on reduced-order and data-driven descriptions of clustering in interacting particle systems~\citep{wehlitz2026data,WehlitzSadeghiMontefuscoetal.2025}.

\new{Through carefully designed numerical experiments, we compare simple scalar observables to illustrate and classify transient clustering dynamics, with particular emphasis on their (non-)monotonicity. 
Whether monotone behaviour is observed depends strongly on the system at hand and on the chosen observable.} For example, if the homogeneous state is asymptotically stable and the initial condition consists of a localised Gaussian peak, should the peak height decrease monotonically in time? 
\new{Peak or maximum density has also
been used in related aggregation models as a simple quantitative
characteristic of aggregate structure
\citep{topaz2006nonlocal,bernoff2020agent,bollinger2016fluids}.
Here, we use peak height as an intuitive geometric indicator of local
concentration and investigate to what extent it reflects the underlying
clustering dynamics. However, peak height is a local geometric quantity
and contains only limited information about the overall system configuration.}
Different density profiles with substantially different energy levels may exhibit the same maximal peak value. Alternatively, many other quantities could serve as order parameters for clustering, including  dispersion measures, Fourier-mode amplitudes, interaction energies, or distances to equilibrium. These observables capture different aspects of clustering, such as the number of clusters, their spatial localisation, separation, or sharpness. While local clustering indicators such as peak height or local mass concentration reflect density variations at specific spatial locations, global quantities such as entropy or interaction energy capture collective effects across the entire domain. Different observables therefore encode complementary information about the transient clustering dynamics and the underlying transport mechanisms.

\new{In this work, we use the standard Wasserstein gradient-flow structure of the
McKean--Vlasov equation and the corresponding decomposition of the free energy
into entropic and interaction contributions as a diagnostic for the transient
dynamics. By tracking the temporal evolution of these contributions, we
distinguish aggregation-dominated, diffusion-dominated, and cooperative
regimes and identify switches between them. 
We compare this energetic classification with geometric clustering observables,
in particular the density peak height, to investigate the non-monotone character
of transient clustering dynamics and the extent to which such observables reflect
the underlying transport mechanisms. In particular, the peak height may increase
even during a diffusion-dominated regime, showing that intuitive geometric
indicators do not necessarily identify the dominant energetic mechanism. This perspective may be relevant not only to
applications in physics or chemistry (see, e.g.,~\citet{taboada2005modeling,rouwhorst2020nonequilibrium}), but also to the
analysis of transient dynamics during crowd evacuation
\citep{bellomo2016human}, the design of control strategies against group
polarisation \citep{zimper2025meanfieldoptimalcontrolstochastic}, or the understanding of generalisation in
machine-learned diffusion models under early stopping
\citep[see, e.g.,][]{li2023generalization}, to mention just a few examples.}

The interacting particle system provides a microscopic description whose collective behaviour is represented at the macroscopic level by the mean-field McKean--Vlasov equation studied in this work. 
We analyse the resulting one-dimensional aggregation--diffusion PDE on the torus with periodic boundary conditions. As representative examples of locally attractive interactions, we consider Morse-type interaction potentials~\mbox{\citep{dorsogna2006, carrillo2019}} as well as the Hegselmann--Krause model~\mbox{\citep{hegselmann2002,
garnier2017consensus,gerber2025formation}}. 
\new{Through numerical experiments, we observe non-monotone evolution of geometric
clustering observables, in particular the density peak height, along trajectories
converging either to homogeneous or to clustered equilibria. To analyse the
underlying transport mechanisms, we identify alternating aggregation- and
diffusion-dominated transient regimes, despite monotone dissipation of the total
free energy. 
These energetic regime changes are only partially reflected in the evolution of the peak height, showing that geometric clustering observables and the analysis of the free-energy contributions provide complementary information about the transient dynamics.}

\paragraph{Outline.} In~\Cref{sec:model_energy}, we introduce the stochastic aggregation--diffusion model together with its mean-field McKean--Vlasov formulation\new{, associated
gradient-flow structure, and the dependence of stationary states on the noise strength.} In~\Cref{sec:transient}, \new{we investigate non-monotone transient clustering dynamics by using the free-energy decomposition to classify aggregation-, diffusion-, and cooperative regimes and present the corresponding numerical experiments.} 
Finally, conclusions and perspectives for future work are discussed in~\Cref{sec:conclusion}.

\section{Model and energetic structure}
\label{sec:model_energy}

We begin in \Cref{sec:model} by introducing the stochastic particle-based aggregation--diffusion system together with its mean-field limit. The associated \new{standard} gradient-flow formulation and free-energy structure are presented in \Cref{sec:energy}. Finally, in \Cref{sec:stationary}, we briefly summarize the dependence of stationary states on the model parameters.

\subsection{Interacting particle system and mean-field limit} \label{sec:model}

We consider a system of $N$ interacting particles moving on the $d$-dimensional torus $\mathbb{T}^d \simeq [-\frac{L}{2},\frac{L}{2})^d$ with periodic boundary conditions. The microscopic dynamics are described by an overdamped Langevin system
\begin{equation}
\label{eq:dX}
\mathrm{d}X_i(t)
= -\frac{1}{N}\sum_{\substack{j=1\\ j\neq i}}^N \nabla U\!\left(X_i(t)-X_j(t)\right)\,\mathrm{d}t
+ \sigma\,\mathrm{d}W_i(t),
\qquad i=1,\dots,N,
\end{equation}
where $X_i(t)\in\mathbb{T}^d$ denotes the position of particle $i$ at time $t$, $U:\mathbb{R}^d\to\mathbb{R}$ is a pairwise interaction potential which is understood to be periodically extended to $\mathbb{T}^d$, $\sigma>0$ is the noise amplitude, and~$(W_i)_{i=1}^N$ are independent $d$-dimensional standard Wiener processes. In contrast to models with an additional external confining potential~\citep{yang2026noisy}, the present system is purely interaction-driven and confined solely through the periodic geometry of the torus, following the setup of \citet{carrillo_long-time_2020}. 

In the mean-field limit $N\to\infty$, the empirical measure
\begin{equation*}
\rho^N(x,t) = \frac{1}{N}\sum_{i=1}^N \delta(x-X_i(t))
\end{equation*}
converges (in an appropriate weak sense and for suitable interaction potentials $U$) to a deterministic particle density $\rho(x,t)$ satisfying a nonlocal nonlinear PDE, also called \textit{McKean--Vlasov PDE},
\begin{equation}
\label{eq:MKV}
\partial_t \rho(x,t)
= \frac{\sigma^2}{2}\,\Delta \rho(x,t) + \nabla\!\cdot\!\bigl(\rho(x,t)\, (\nabla U * \rho)(x,t)\bigr),
\end{equation}
where $*$ denotes convolution on $\mathbb{T}^d$~\citep{dawson1983critical, gartner1988mckean, sznitman1991}. The first term models diffusion due to thermal fluctuations, while the second term represents nonlocal interactions between particles. We analyse the clustering behaviour of the interacting particle system by means of its mean-field limit. 

\paragraph{Specialization to one dimension and choice of interaction.}
In the remainder of this work, we restrict attention to the one-dimensional case $d=1$, for which all analytical considerations and numerical experiments are performed. As a representative example of an interaction potential, we consider a Morse-type potential of the form
\begin{equation}
\label{eq:Morse}
U(x) = -C_a\,\mathrm{exp}\left(-|x|/\ell_a\right) + C_r\,\mathrm{exp}\left(-|x|/\ell_r\right),
\end{equation}
with attraction and repulsion strengths $C_a,C_r>0$ and characteristic length scales $\ell_a,\ell_r>0$.  

In many applications, the Morse potential is used to model short-range repulsion and long-range attraction in swarming phenomena \citep{dorsogna2006, leverentz2009, carrillo2019}. Here, however, we select parameter values such that the potential becomes locally attractive, leading to short-range interactions that promote cluster formation. Moreover, we employ a periodized version of the Morse potential\new{~\eqref{eq:Morse} on the torus $\mathbb{T}\simeq[-\frac{L}{2}, \frac{L}{2})$ (see, e.g., \cite{bandegi2018}), given by
\begin{equation}
\label{eq:U_per_Morse}
    U_{\text{per}}(x) = - C_a \, \frac{\cosh\left(\frac{|x| - \frac{L}{2}} {l_a}\right)}{\sinh\left(\frac{L}{2l_a}\right)} \, + \,  C_r \, \frac{\cosh \left(\frac{|x| - \frac{L}{2}} {l_r}\right)}{\sinh\left(\frac{L}{2l_r}\right)} \, .
\end{equation}
This periodized version is obtained by summing the original potential over all periodic images, i.e., over positions shifted by integer multiples of the period $L$, and evaluating the resulting sums using the geometric series.}

\new{Because the Morse potential~\eqref{eq:Morse} depends on $|x|$, it is in
general not differentiable at the origin. More precisely, unless
$C_a/\ell_a=C_r/\ell_r$, the one-sided limits of $\nabla U$ at $x=0$
differ, so that the interaction force $-\nabla U$ has a finite jump there.
We define $\nabla U(0):=0$, thereby obtaining an everywhere-defined bounded
force, which nevertheless remains discontinuous. 
Consequently, the classical propagation-of-chaos result of
\citet{sznitman1991}, which assumes Lipschitz-continuous interaction forces,
does not apply directly. 
The more recent result of~\citet{jabin2018}, however, allows interaction
forces of much lower regularity, formulated in terms of weak distributional
bounds, and in particular includes bounded forces with jump discontinuities
such as the Morse force considered here. 
Related propagation-of-chaos results for singular interaction kernels,
including uniform-in-time estimates, have subsequently been obtained
by~\citet{guillin2025}.}

While the numerical results presented below focus primarily on this specific interaction kernel, the qualitative phenomena discussed---including phase coexistence and alternating transient regimes---are not specific to the Morse potential 
and extend to a broader class of locally attractive interactions. 
To further illustrate this robustness, we additionally consider the Hegselmann--Krause model~\citep{hegselmann2002, garnier2017consensus, gerber2025formation} in Appendix~\ref{sec:Hegselmann--Krause model}.

\subsection{Free energy and gradient-flow formulation} \label{sec:energy}

\new{We briefly recall the standard Wasserstein gradient-flow formulation of
the PDE~\eqref{eq:MKV} and its associated thermodynamic structure~\citep{jordan1998variational,ambrosio2008gradient,santambrogio2015}. The
corresponding free-energy functional}
\begin{equation}
\label{eq:free_energy}
\mathcal{F}[\rho]
=  \underbrace{\frac{\sigma^2}{2}\int_{\mathbb{T}^d} \rho(x)\log\rho(x)\,\mathrm{d}x}_{:=\mathcal{F}_{\text{ent}}[\rho]} \, + \, \underbrace{\frac12 \int_{\mathbb{T}^d}\!\!\int_{\mathbb{T}^d}
\rho(x)\rho(y)\,U(x-y)\,\mathrm{d}x\,\mathrm{d}y}_{:=\mathcal{F}_{\text{int}}[\rho]}
\end{equation}
\new{is} given by the sum of an entropic contribution~$\mathcal{F}_{\text{ent}}[\rho]$ and a nonlocal interaction-energy contribution~$\mathcal{F}_{\text{int}}[\rho]$.
More concretely, we have 
\begin{equation}
\label{eq:gradient_flow}
    \partial_t \rho
    =
    \nabla \cdot \left( \rho \nabla \mu \right),
\end{equation}
where
\begin{equation}
\label{eq:mu}
\mu
=
\frac{\delta\mathcal F}{\delta\rho}
=
\frac{\sigma^2}{2}(1+\log\rho)
+
U*\rho .   
\end{equation}
denotes the first variation (variational derivative) of the free-energy functional $\mathcal F$ and will be referred to as the \textit{chemical potential}; it represents the (Helmholtz) free energy cost associated with adding an infinitesimal amount of particles to the system. \new{Note that in the chemical potential~\eqref{eq:mu}, the additive constant $\frac{\sigma^2}{2}$ 
does not affect the dynamics, since only $\nabla \mu$ enters~\eqref{eq:gradient_flow}.}

The chemical potential determines the direction of transport through the force field $-\nabla\mu$, while the associated mass flux is given by
\[
J=-\rho\nabla\mu.
\]
Importantly, the density enters only at the level of the flux. Regions of large density therefore contribute more strongly to transport even if the gradient of $\mu$ is moderate. This distinction between chemical potential, its gradient, and the density-weighted flux will play a key role in the interpretation of the transient dynamics below.

\paragraph{Energy dissipation.} Along sufficiently regular solutions of~\eqref{eq:MKV}, the free energy dissipates according to
\begin{align}
     \frac{d}{dt}\mathcal{F}[\rho(t)]
     =
     -\int_{\mathbb{T}^d}
     \rho(x,t)\left|\nabla \mu(x,t)\right|^2
     \,dx      =
     -\int_{\mathbb{T}^d}
     \rho \,
     \left|
        \frac{\sigma^2}{2}\nabla\log\rho \, + \,  \nabla U * \rho
     \right|^2
     \,dx
     \le 0. \label{eq:energy_dissipation}
\end{align}
Thus, the total free energy decreases monotonically along solutions of the PDE, and stationary states correspond to critical points of $\mathcal F$. 
In contrast to a formal $L^2$-gradient flow, the Wasserstein geometry naturally incorporates mass conservation and transport. As a consequence, the dynamics evolve through redistribution of mass rather than through pointwise creation or destruction.

The competition between attractive interactions encoded in the potential $U$ and entropic diffusion scaled by $\sigma$ determines both the transient and long-time behaviour of the system. The decomposition of the free energy into entropic and interaction contributions is reminiscent of the Jordan--Kinderlehrer--Otto (JKO) variational formulation of reversible Fokker--Planck equations~\citep{jordan1998variational,ambrosio2008gradient}. In the present McKean--Vlasov setting, the potential contribution is replaced by a nonlocal, density-dependent interaction energy, which causes the thermodynamic landscape itself to evolve with the density. 

Although the monotone decrease of the total free energy is independent of the parameter values, the relative influence of the interaction energy $\mathcal{F}_{\text{int}}[\rho]$ and the entropy $\mathcal{F}_{\text{ent}}[\rho]$ may vary significantly along the transient dynamics. Understanding the competition between these two contributions forms the central focus of the analysis in~\Cref{sec:transient}.

We briefly mention that a similar interplay between entropy and \new{internal energy} contributions during the relaxation to equilibrium has been observed in non-reversible, hypocoercive Fokker--Planck equations associated with linear stochastic differential equations of Ornstein--Uhlenbeck type \citep{arnold2014sharp,neureither2017time}. Here, the interaction potential $U$ is zero, and the free energy~$\mathcal{F}(\rho)$ is given by the relative entropy $\mathcal{F}(\rho)=D(\rho|\nu)$ between $\rho$ and the Gibbs-Boltzmann density~$\nu=\exp(-2V/\sigma^2)$ with $V$ being a (quadratic) confining potential. Interestingly, the competition between 
different relaxation mechanisms leads to a monotone, but not strictly monotone decay of the free energy, which is not present in reversible Fokker--Planck equations.

\subsection{Stationary states and their dependence on the noise strength} 
\label{sec:stationary}

The stability properties of stationary states are determined by the variational structure of the free-energy functional $\mathcal F$. In particular, local minimizers of $\mathcal{F}$ are dynamically stable, while saddle points or maxima are unstable. The qualitative behaviour of the system is therefore governed by the convexity properties of the free energy, which depend on the interplay between the nonlocal interaction energy and the entropic diffusion term.

For interaction kernels whose Fourier representation contains negative modes, the free energy may become nonconvex as the noise strength $\sigma$ is decreased~\citep{carrillo_long-time_2020}. On periodic domains, such negative Fourier modes favour spatially non-uniform perturbations of the homogeneous state and may lead 
to the emergence of clustered stationary solutions. As a consequence, two critical noise strengths
\[
\sigma_\#, \, \sigma_c>0\quad\text{with}\quad\sigma_\# <  \sigma_c
\]
typically arise. 
For $\sigma > \sigma_c$, the spatially homogeneous density is the unique stable equilibrium. 
As $\sigma$ decreases below $\sigma_c$, additional non-uniform 
stationary states corresponding to clustered particle configurations emerge, leading to a coexistence regime in which both homogeneous and clustered equilibria are locally stable. Finally, for $\sigma < \sigma_\#$, the homogeneous equilibrium loses stability. 

In the coexistence regime, the asymptotic state depends sensitively on the initial condition. Such bistability and loss of stability of the homogeneous state are well known features of McKean--Vlasov dynamics on periodic domains \citep{carrillo_long-time_2020,chayes2010mckean}. This equilibrium structure provides the theoretical background for the numerical investigations presented in the following section, where we focus on the transient dynamics and, in particular, on the occurrence of alternating aggregation- and diffusion-dominated regimes for non-uniform initial conditions.

\section{\new{Non-monotone transient clustering dynamics}}
\label{sec:transient}

\new{We now turn to the non-monotone transient dynamics of the
aggregation--diffusion system. In particular, natural geometric clustering
observables such as the density peak height may evolve non-monotonically along
trajectories approaching equilibrium.}

\new{To understand the mechanisms underlying this behaviour, we use the standard decomposition of the
free energy into interaction and entropic contributions to distinguish
aggregation-dominated, diffusion-dominated, and cooperative regimes. \Cref{sec:energetic} introduces this energetic classification, while \Cref{subsec:Numerical investigation}
compares it with the evolution of the density peak height in representative
numerical examples.}

\subsection{\new{Energetic classification of transient regimes}}
\label{sec:energetic}

\new{For the numerical analysis below, we classify the transient dynamics according
to the signs of the temporal changes of the entropic and interaction
contributions to the free energy.}
The total free energy satisfies
\begin{equation}
\frac{d}{dt}\mathcal F
= \frac{d}{dt}\mathcal F_{\text{ent}} + \frac{d}{dt}\mathcal F_{\text{int}}
\le 0,
\end{equation}
with equality only at equilibrium; see~\Cref{sec:energy}. The competition between aggregation and diffusion is therefore encoded in the relative behaviour of the interaction and entropic contributions.

\begin{itemize}
    \item \new{The system is in an \emph{aggregation-dominated} regime when} 
\begin{equation}
  \frac{d}{dt} \mathcal{F}_{\text{ent}} >0, \qquad \frac{d}{dt} \mathcal{F}_{\text{int}} <0,
\end{equation}
corresponding to interaction-driven mass concentration despite an entropic cost.
\item \new{The system is in a \emph{diffusion-dominated} regime when} 
\begin{equation}
       \frac{d}{dt} \mathcal{F}_{\text{ent}} <0, \qquad \frac{d}{dt} \mathcal{F}_{\text{int}} >0, 
\end{equation}
corresponding to entropically driven spreading and an increase in interaction energy.
\item \new{The system is in a \emph{cooperative} relaxation regime when}
\begin{equation}
       \frac{d}{dt} \mathcal{F}_{\text{ent}} <0, \qquad \frac{d}{dt} \mathcal{F}_{\text{int}} <0,
\end{equation}
\new{so that both contributions decrease simultaneously.}
\end{itemize}
\new{\emph{Alternating transient regimes} arise when the dynamics switch between these different sign structures over time, even though the total free energy decreases monotonically throughout the evolution.}
Depending on the initial condition and the noise strength, such trajectories may converge either to the homogeneous equilibrium or to a clustered stationary state. Two representative examples are presented in~\Cref{sec:alternating}: in the first case (\Cref{ex1}), the system eventually relaxes to the homogeneous equilibrium, whereas in the second case (\Cref{ex:twoGaussians}) it converges to a stable clustered state. 

Before turning to these alternating regimes, we first consider in~\Cref{sec:single-dominance} the simpler single-dominance regimes, where the dynamics remain either diffusion-dominated or aggregation-dominated throughout the entire time evolution.

\subsection{Numerical investigation}
\label{subsec:Numerical investigation}

As a representative example of an interaction kernel, we consider the Morse-type potential~\eqref{eq:Morse} on the one-dimensional torus $\mathbb{T}\simeq[-\frac{L}{2}, \frac{L}{2})$. Throughout our numerical experiments, we use the parameter values
\begin{align}
\label{eq:parameter_values}
L=5,\qquad
C_a=4,\qquad
C_r=1,\qquad
\ell_a=0.025\,L,\qquad
\ell_r=0.01\,L.
\end{align}
This choice yields a locally attractive interaction force with a bounded jump discontinuity at the origin; \new{there is no rigorous cut-off of the potential beyond a certain radius, but the potential decays rapidly to zero and hence exerts no significant attractive force at larger distances.}

Moreover, the corresponding interaction kernel possesses negative Fourier modes on the periodic domain, which according to the discussion in~\Cref{sec:stationary} allows for loss of stability of the homogeneous equilibrium and the emergence of clustered stationary states~\citep{carrillo_long-time_2020}.

For these parameters, the qualitative behaviour of the system depends strongly on the noise strength~$\sigma$. In particular, we distinguish the following regimes:
\begin{itemize}
    \item $\sigma<\sigma_{\#}\approx0.59$: the homogeneous equilibrium is unstable and the dynamics converge towards a clustered stationary state. The value of $\sigma_{\#}$ is computed from~\cite[Eq.~(3.5)]{carrillo_long-time_2020}.
    \item $\sigma_{\#}<\sigma<\sigma_c\approx0.86$: coexistence regime in which both the homogeneous and a clustered equilibrium are locally stable. \new{The critical value $\sigma_c \approx 0.86$ was estimated numerically from long-time PDE simulations; details of the procedure are given in Appendix~\ref{app:sigmac}.}
    \item $\sigma>\sigma_c\approx0.86$: the homogeneous equilibrium is the unique stable stationary state.
\end{itemize}

\paragraph{Numerical solution of the aggregation--diffusion equation} We solve the McKean--Vlasov PDE~\eqref{eq:MKV} using the positivity-preserving and energy-dissipating finite-volume scheme introduced in~\citet{carrillo2015,bailo2020}. Throughout all numerical experiments, we choose a time step of~$dt = 0.001$ and discretize the one-dimensional torus $\mathbb{T} \simeq \left[-\frac{L}{2}, \frac{L}{2}\right)$ using a uniform grid with $2^9$ cells.

\subsubsection{Single-dominance regimes} 
\label{sec:single-dominance}

We first consider parameter regimes in which, for suitable initial conditions, the dynamics evolve without alternating aggregation- and diffusion-dominated phases. These cases serve as reference scenarios for the alternating transient regimes discussed below.

\begin{figure}[t]
 \begin{subfigure}{0.41\textwidth}
        \centering
    \includegraphics[width=\textwidth]{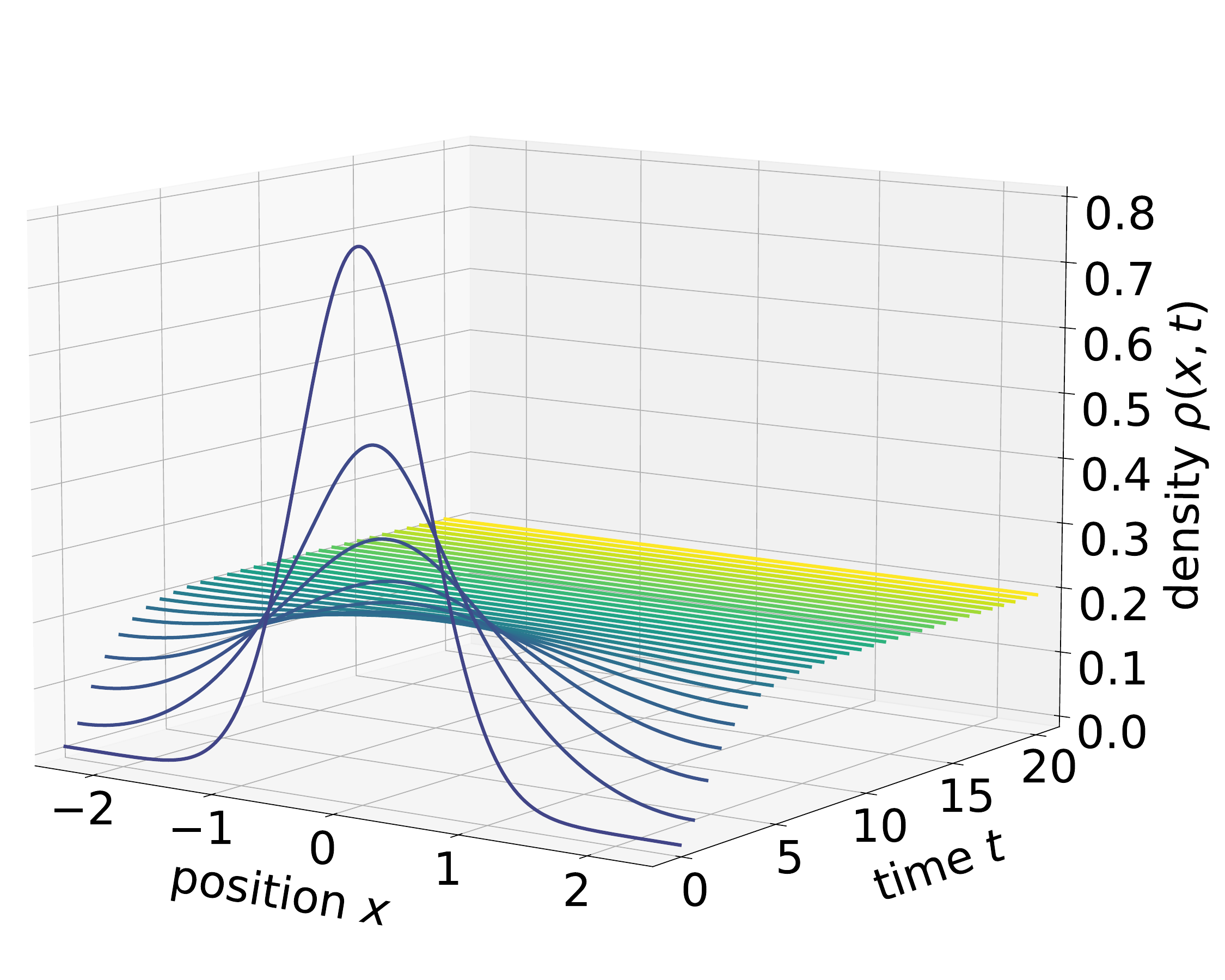}
     \caption{Density evolution $\rho(x,t)$}
     \label{fig:oneGaussian_diffusion_a}
    \end{subfigure}
    \begin{subfigure}{0.58\textwidth}
        \centering
    \includegraphics[width=\textwidth]{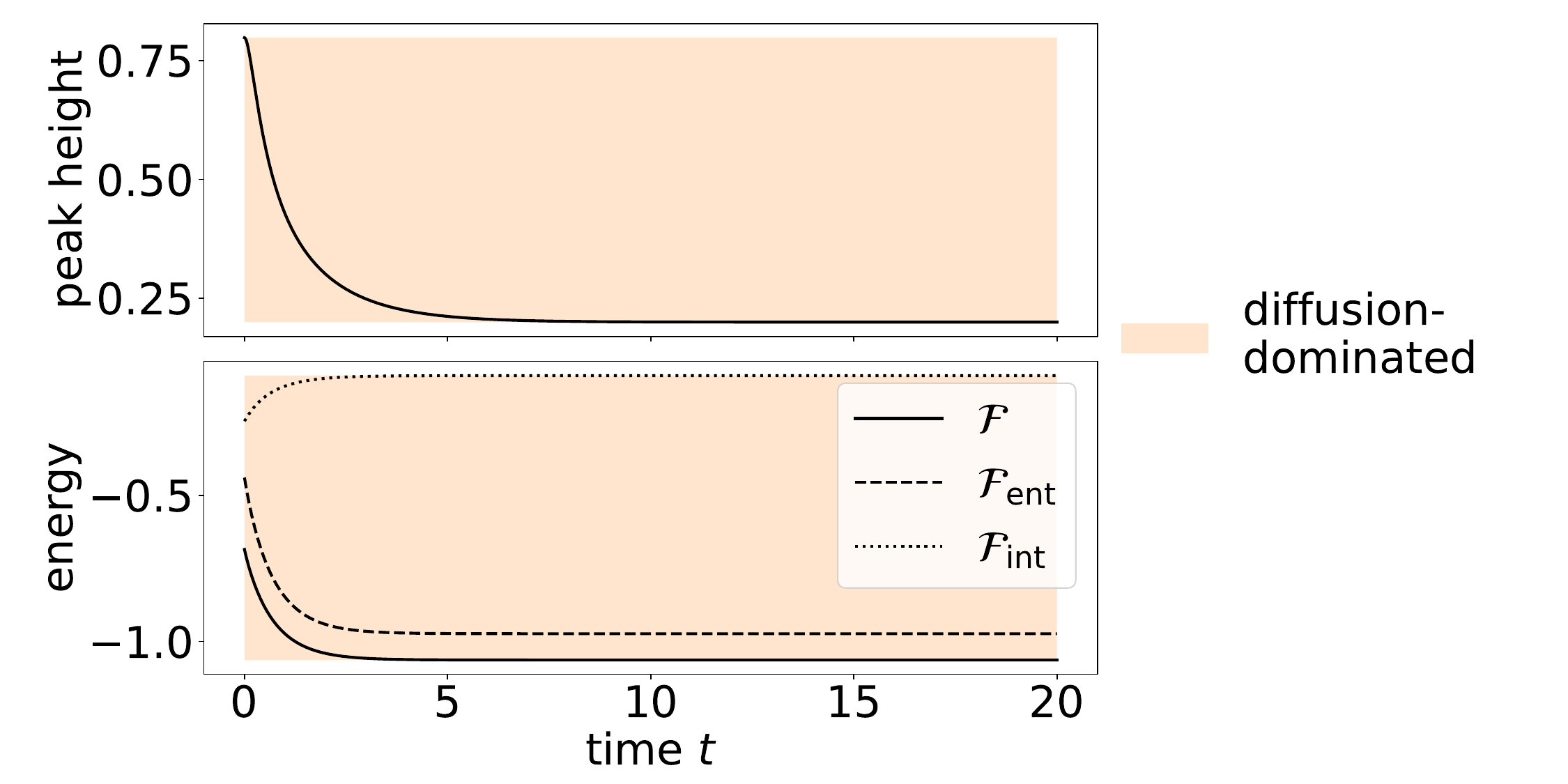}
     \caption{Peak height and free-energy contributions}
     \label{fig:oneGaussian_diffusion_b}
    \end{subfigure}
    \caption{\textbf{Pure diffusion-dominated relaxation to the homogeneous equilibrium.}
    Numerical solution of the PDE~\eqref{eq:MKV} for $\sigma=1.1$ starting from
    $\rho_0\sim \mathcal{N}_{\mathrm{per}}(0,0.5^2)$. 
    (a) Time evolution of the spatial density profiles $\rho(x,t)$, showing spreading towards the homogeneous equilibrium.
    (b) Time evolution of the peak height and free-energy contributions. The peak height decreases monotonically, while the free-energy contributions indicate diffusion-dominated dynamics throughout the evolution.
    }
    \label{fig:oneGaussian_diffusion}
    \end{figure} 

For sufficiently large values of the noise strength~$\sigma$, the evolution is diffusion-dominated. In this regime, the entropic contribution to the free energy decreases while the interaction energy increases, and the density relaxes smoothly towards the homogeneous equilibrium. This behaviour is illustrated in Figure~\ref{fig:oneGaussian_diffusion} for a locally attractive Morse-type interaction potential with parameter values given in~\eqref{eq:parameter_values}. The initial state $\rho_0 = \rho(\cdot,0)$ is chosen as a periodized Gaussian with mean $0$ and variance $0.5^2$,
\[\rho_0\sim \mathcal{N}_{\text{per}}(0,0.5^2).\]
The noise strength $\sigma=1.1>\sigma_c\approx0.86$ ensures that the homogeneous equilibrium is the unique stable stationary state. 
\new{\Cref{fig:oneGaussian_diffusion_a} shows the density evolution in a waterfall
representation, where each curve corresponds to the density profile
$\rho(\cdot,t)$ at a given time. \Cref{fig:oneGaussian_diffusion_b} shows the corresponding peak height
(upper panel) and free-energy contributions (lower panel) over time.}

\begin{figure}
 \begin{subfigure}{0.41\textwidth}
        \centering
    \includegraphics[width=\textwidth]{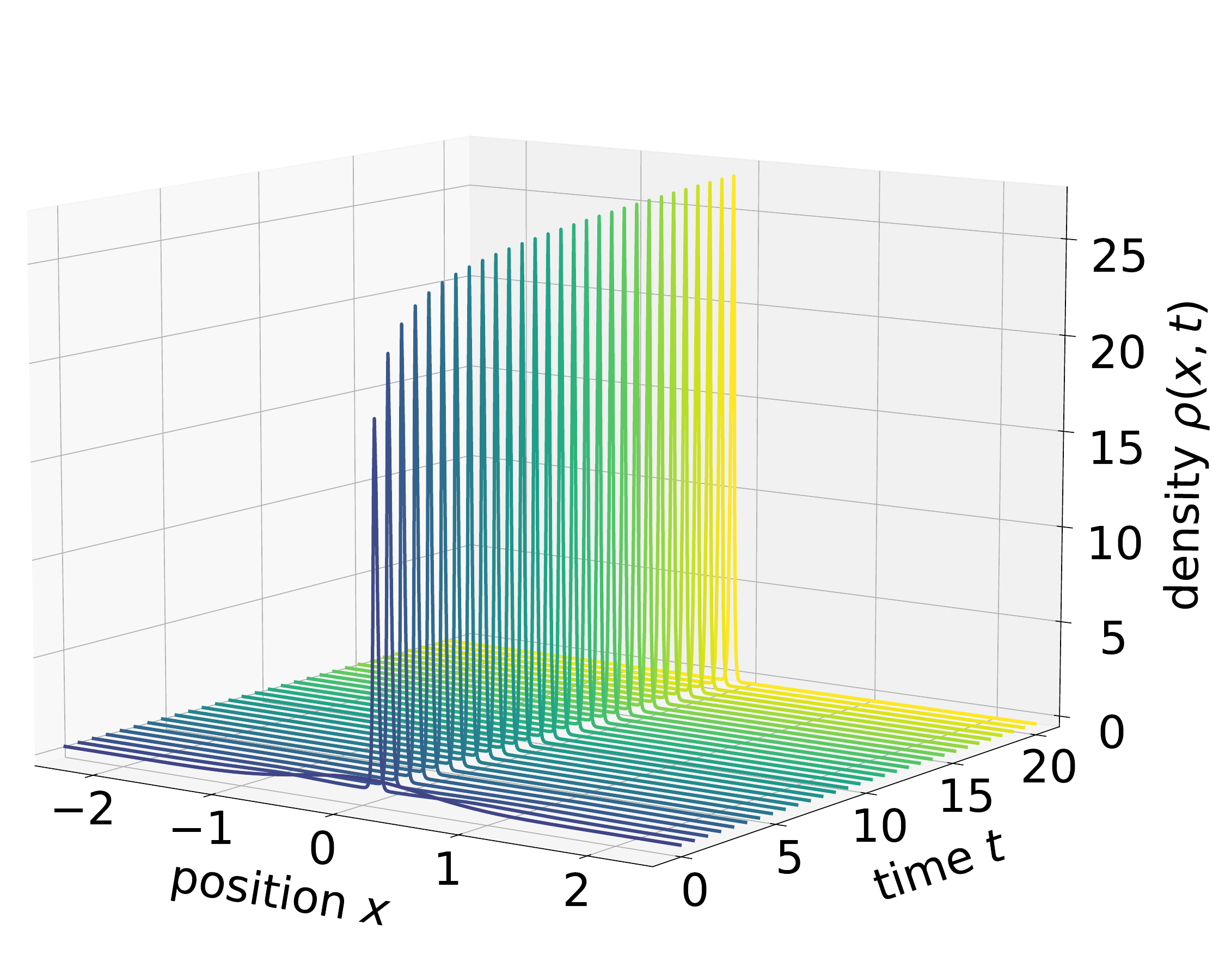}
     \caption{Density evolution $\rho(x,t)$}
    \end{subfigure}
    \begin{subfigure}{0.58\textwidth}
        \centering
    \includegraphics[width=\textwidth]{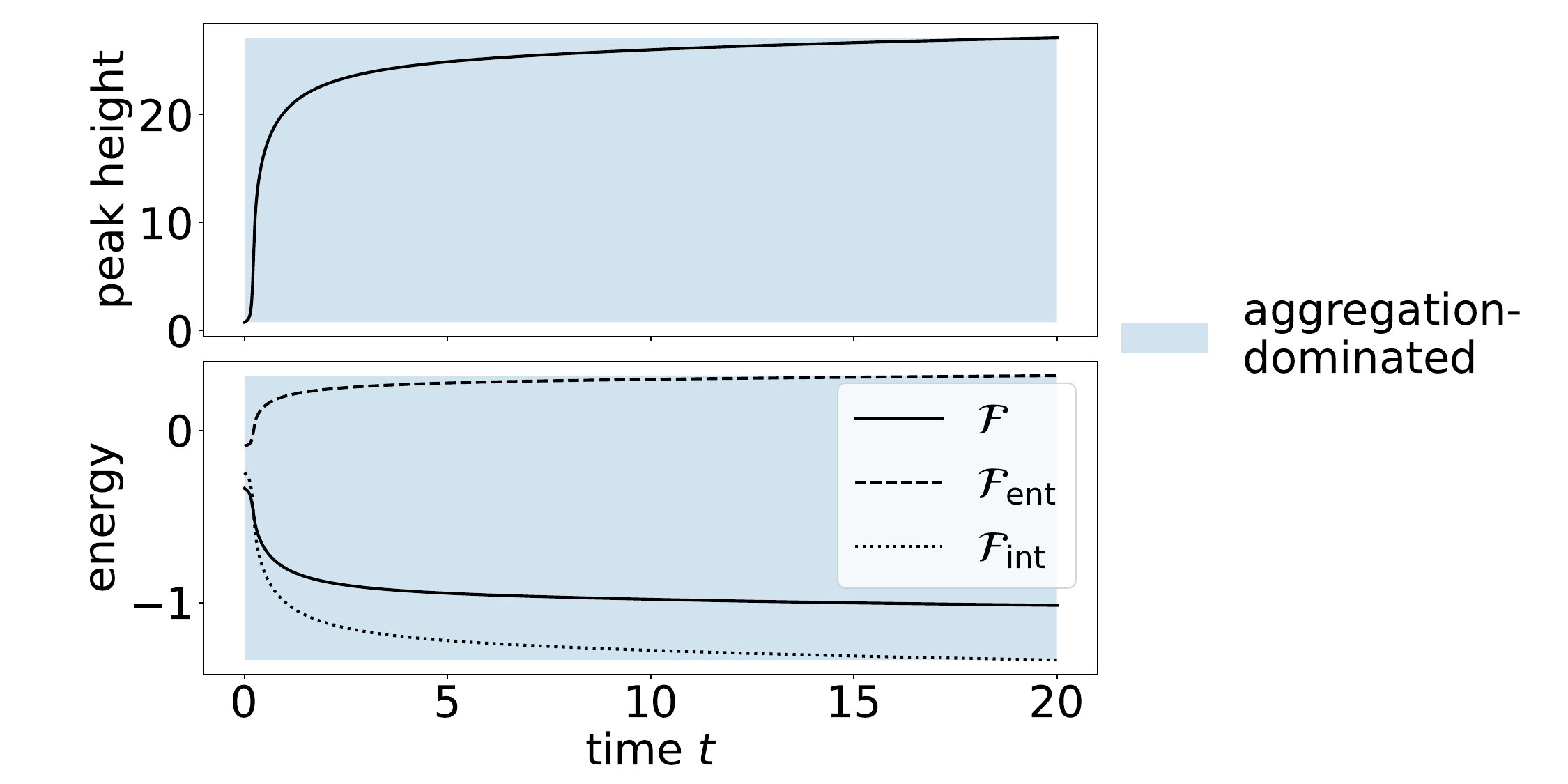}
     \caption{Peak height and free-energy contributions}
    \end{subfigure}
    \caption{\textbf{Pure aggregation-dominated relaxation to the clustered equilibrium.} Numerical solution of the PDE~\eqref{eq:MKV} for $\sigma= 0.5$ starting from~$\rho_0\sim \mathcal{N}_{\text{per}}(0,0.5^2)$. 
    (a) Time evolution of the spatial density profiles $\rho(x,t)$, showing concentration towards a clustered equilibrium.
    (b) Time evolution of the peak height and free-energy contributions. The peak height increases monotonically, while the free-energy contributions indicate aggregation-dominated dynamics throughout the evolution.
    }
    \label{fig:oneGaussian_aggregation}
    \end{figure}

Conversely, for sufficiently small values of $\sigma$, aggregation dominates the dynamics. In this case, the interaction energy decreases while the entropy increases, and the density evolves directly towards a clustered equilibrium. \Cref{fig:oneGaussian_aggregation} shows an example with $\sigma=0.5<\sigma_{\#}\approx0.59$, for which the homogeneous equilibrium is unstable and the dynamics converge towards a stable single-peak clustered state.

These single-dominance regimes provide a baseline for the intermediate noise regime, where alternating aggregation- and diffusion-dominated phases emerge.

\subsubsection{Alternating transient regimes} \label{sec:alternating}

Alternating transient regimes arise when the dynamics switch repeatedly between aggregation- and diffusion-dominated behaviour.
Depending on the initial condition and noise strength, such alternating transient dynamics may relax either to the homogeneous equilibrium or to a clustered stationary state. We discuss these two scenarios separately below. In both cases, the Morse-type interaction potential~\eqref{eq:Morse} with parameter values~\eqref{eq:parameter_values} is considered.

\FloatBarrier 

\begin{ex}[Alternating transient regimes relaxing to the homogeneous equilibrium] \label{ex1} As before, we start from the initial condition
\[
\rho_0=\rho(\cdot,0)
\]
that is chosen as a periodized Gaussian with mean $0$ and variance $0.5^2$. 
The noise strength
\[
\sigma=0.838
\]
lies in the coexistence regime close to the critical value $\sigma_c\approx0.86$, where both the homogeneous equilibrium and a clustered stationary state are locally stable. For the present initial condition, the dynamics converge to the homogeneous equilibrium. However, for different initial states, convergence towards the clustered stationary state can also be observed; see~\Cref{fig:oneGaussian_start0.4} in the Appendix. 

\begin{figure}
\begin{subfigure}{0.41\textwidth}
    \centering
\includegraphics[width=\textwidth]{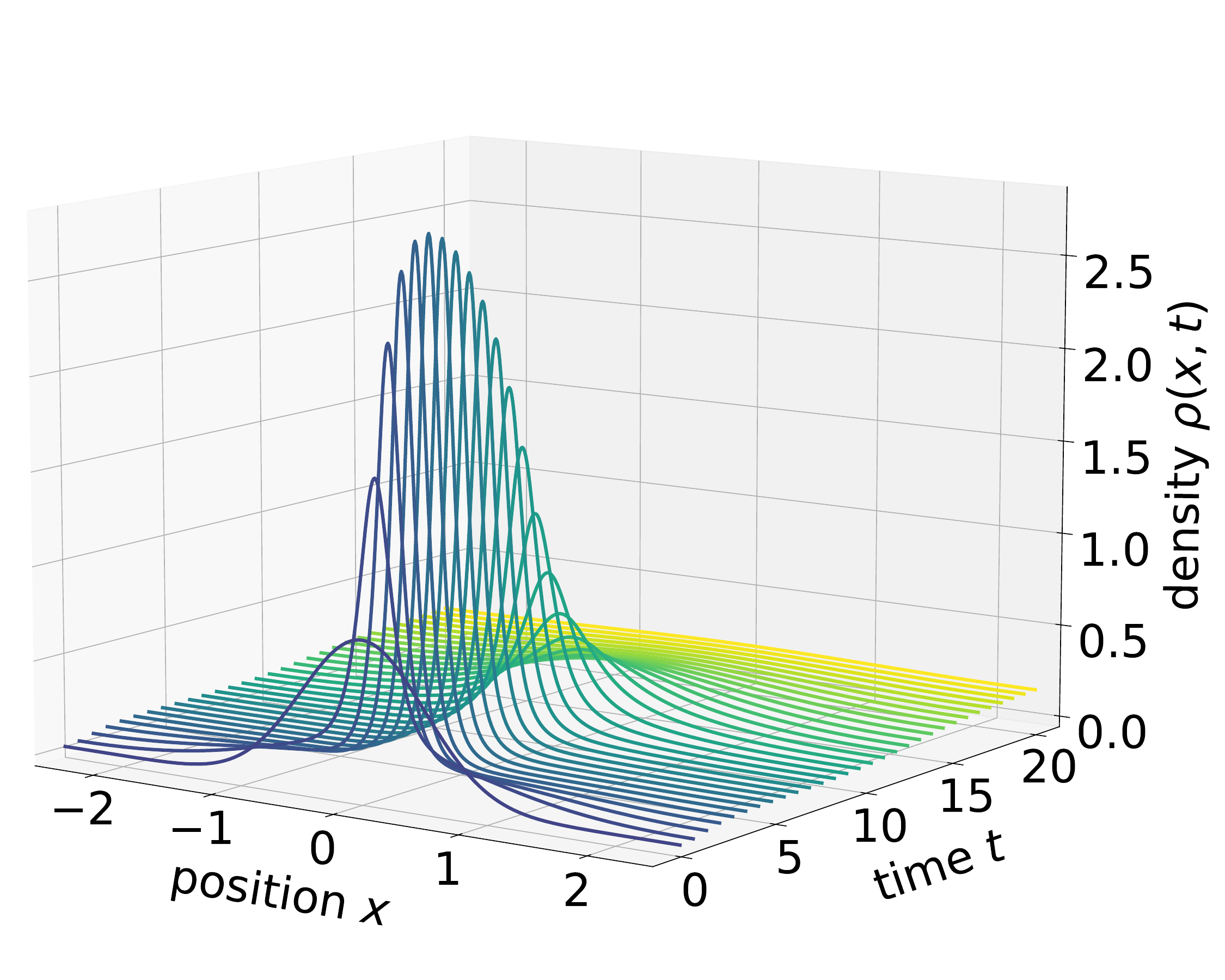}
 \caption{Density evolution $\rho(x,t)$}\label{fig:oneGaussian_rho_3d}
\end{subfigure}
         \begin{subfigure}{0.58\textwidth}
    \centering
\includegraphics[width=\textwidth]{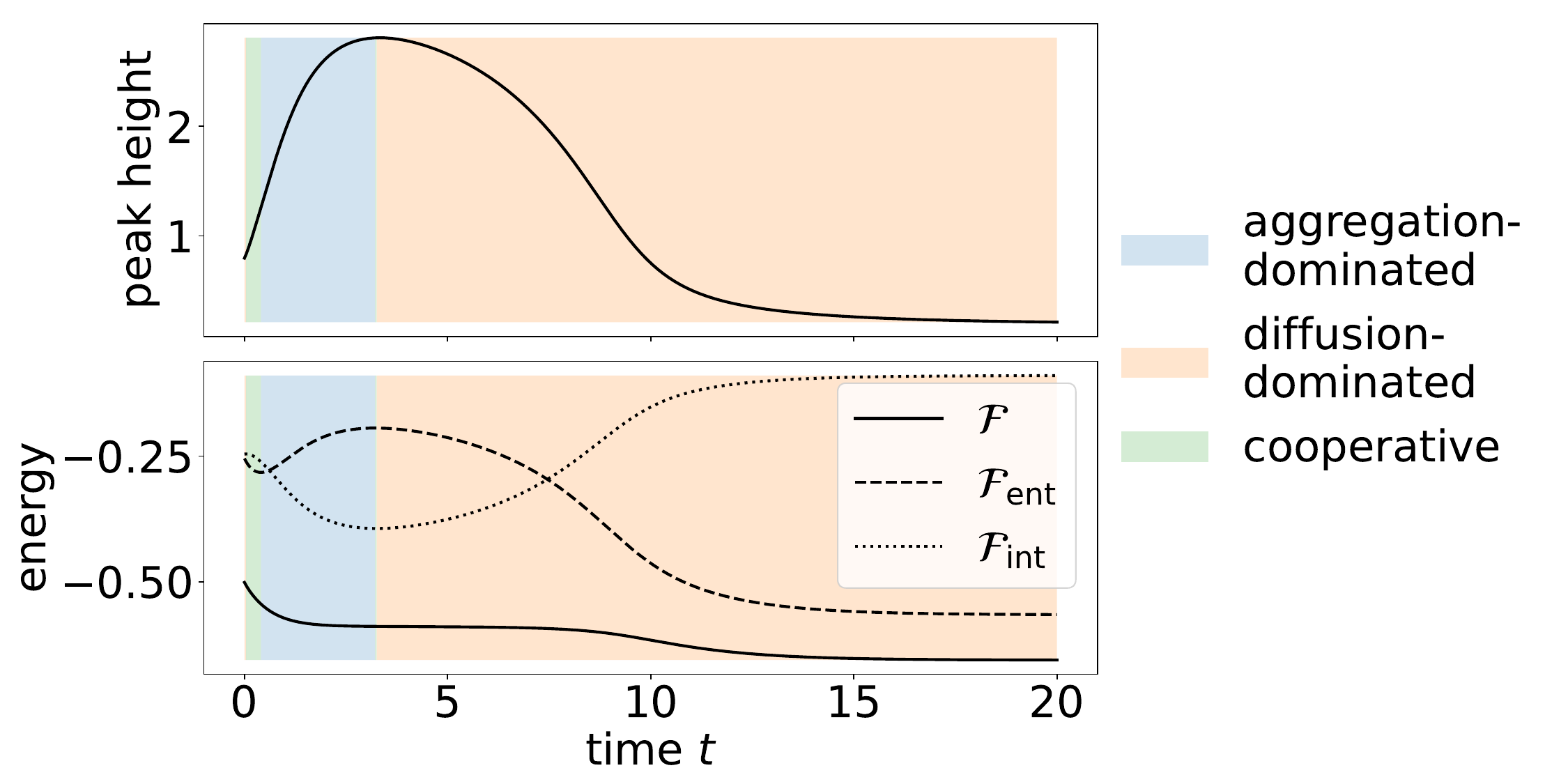}
 \caption{Peak height and free-energy contributions}\label{fig:oneGaussian_F_peak}
\end{subfigure}
\caption{\textbf{Alternating transient regimes during relaxation to the homogeneous equilibrium.} Numerical solution of the PDE~\eqref{eq:MKV} for $\sigma= 0.838$ starting from~$\rho_0\sim \mathcal{N}_{\text{per}}(0,0.5^2)$, see~\Cref{ex1}. 
(a) Time evolution of the spatial density profiles $\rho(x,t)$, showing initial concentration followed by spreading towards the homogeneous equilibrium.
(b) Time evolution of the peak height and free-energy contributions. The peak height first increases and subsequently decreases, while the free-energy contributions reveal a transition from aggregation-dominated to diffusion-dominated dynamics. }
\label{fig:oneGaussian}
\end{figure} 

\new{\Cref{fig:oneGaussian} shows the temporal evolution of the density $\rho$, of the peak height, and of the free energy contributions.} Starting from the non-uniform initial condition, the density first develops a pronounced cluster before spreading towards the homogeneous equilibrium (\Cref{fig:oneGaussian_rho_3d}). Accordingly, the peak height shown in the upper
panel of \Cref{fig:oneGaussian_F_peak} initially increases and later decreases towards the uniform value.

A more detailed picture emerges from the evolution of the free-energy contributions shown in the lower panel of \Cref{fig:oneGaussian_F_peak}. 
\new{Initially, the dynamics are diffusion-dominated ($t\lesssim0.04$), before entering a cooperative relaxation regime in which both entropy and interaction energy decrease ($0.04\lesssim t\lesssim0.4$).} 
\sout{After a short initial diffusion-dominated phase ($t\lesssim0.04$), the dynamics enter a cooperative relaxation regime in which both entropy 
and interaction energy decrease ($0.04\lesssim t\lesssim0.4$).} This is followed by an 
aggregation-dominated phase characterised by decreasing interaction energy and increasing entropy ($0.4\lesssim t\lesssim3.2$). At later times ($t\gtrsim3.2$), the trend reverses and the dynamics become diffusion-dominated again, leading to relaxation towards the homogeneous equilibrium. 
The transitions between these regimes are accompanied by pronounced plateaus in the total free energy, indicating extended time intervals during which the competing energetic contributions nearly balance each other. \hfill$\diamond$
\end{ex}

\begin{figure}
    \centering
    \includegraphics[width=\textwidth]{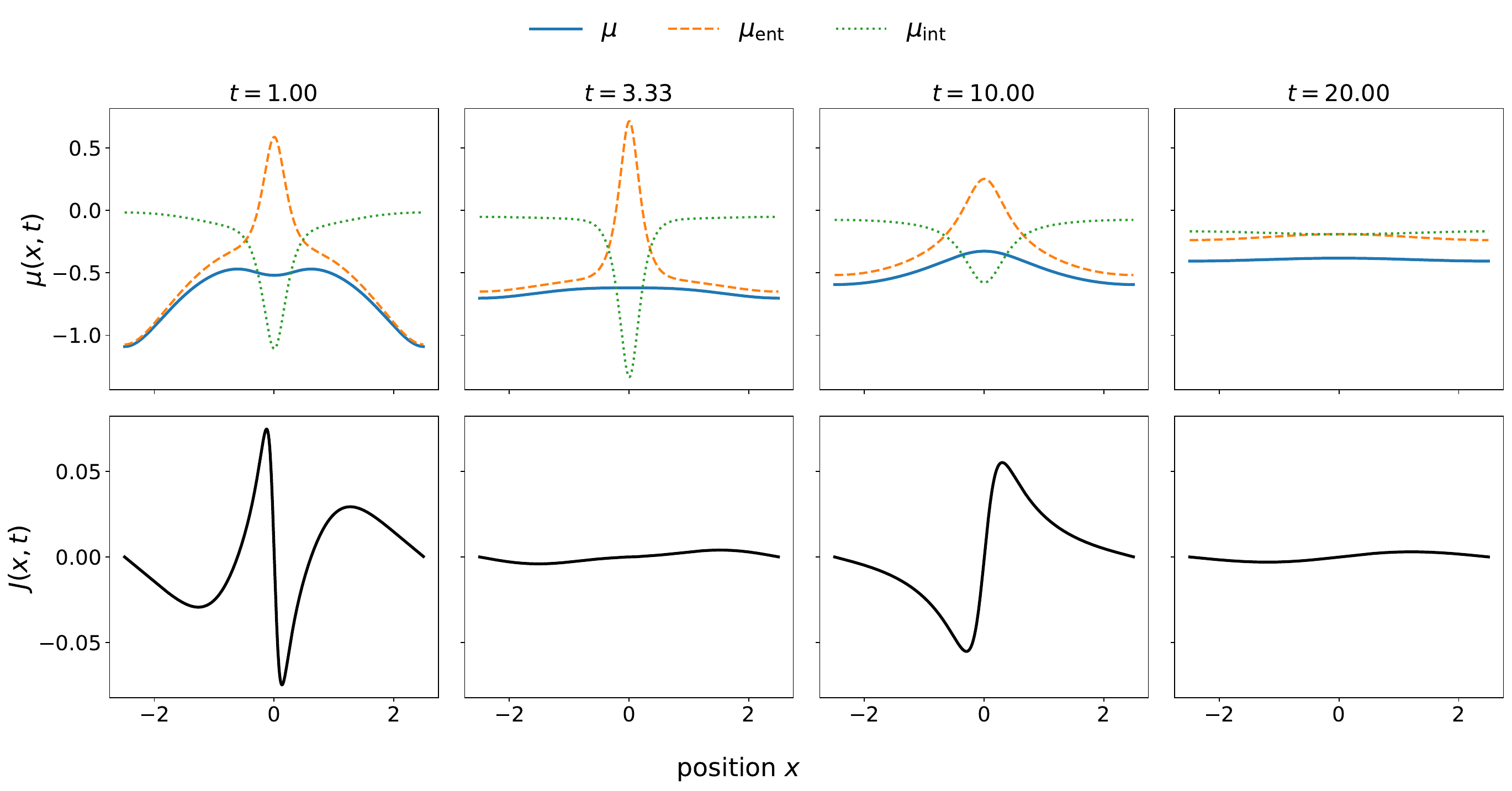}
\caption{\textbf{Local transport mechanisms during the non-monotone transient of~\Cref{ex1}.}
Chemical potential $\mu=\frac{\delta \mathcal{F}}{\delta \rho}$ (top row)
and associated flux $J=-\rho\nabla\mu$ (bottom row) at four representative
times. The chemical potential is decomposed as
$\mu=\mu_{\rm ent}+\mu_{\rm int}$.
At $t=1$, a local minimum of $\mu$ near the centre induces inward transport;
at the peak-height maximum, $t=3.33$, $\mu$ is nearly flat near the centre 
and the flux is small; at $t=10$, during relaxation
towards the homogeneous equilibrium, a local maximum induces outward transport;
and at $t=20$, close to equilibrium, $\mu$ is nearly spatially constant and
the flux is correspondingly small.} 
\label{fig:oneGaussian_2}
\end{figure}

\new{\paragraph{Thermodynamic interpretation of the peak-height dynamics.}
\new{The change from increasing to decreasing peak height can be understood from the local transport induced by the chemical potential. To illustrate
this mechanism, we examine the chemical potential and the associated flux at representative times: \Cref{fig:oneGaussian_2} shows the chemical potential $\mu$ and the associated flux
\[
J=-\rho\nabla\mu
\]
at four representative times during the non-monotone transient in~\Cref{ex1}.
At $t=1$, the interaction contribution creates a local minimum of $\mu$ near
the center. Since the flux is directed down the chemical-potential gradient
and weighted by $\rho$, this produces inward transport in the high-density
central region and hence a further increase of the peak height.
At $t=3.33$, when the peak height is maximal, the spatial variations of the
entropic and interaction contributions largely compensate near the center,
leaving the total chemical potential almost flat and the net flux small.
At $t=10$, the centre has become a local maximum of $\mu$, so that the flux
is directed outward and the peak decreases as the solution relaxes towards
the homogeneous equilibrium. By $t=20$, the system is close to equilibrium, in the sense that 
the chemical potential is nearly spatially constant and the remaining flux is
small.
Thus, the change from increasing to decreasing peak height is accompanied by
a reversal of the local chemical-potential landscape and the associated mass
flux, followed by their decay as equilibrium is approached.}}

\quad 

\new{The relation between peak-height evolution and energetic dominance is,
however, not one-to-one. As shown in \Cref{fig:oneGaussian_start0.6} in Appendix~\ref{app:initial_state}, even a
trajectory that is diffusion-dominated in the energetic sense can exhibit
an initial increase of the peak height before the density subsequently
spreads out. A comparison with the second moment and with a two-dimensional description combining second moment and peak height is provided in~Appendix~\ref{sec:choice_and_number_of_observables}.}

\begin{figure}
    \centering
    \includegraphics[width=\textwidth]{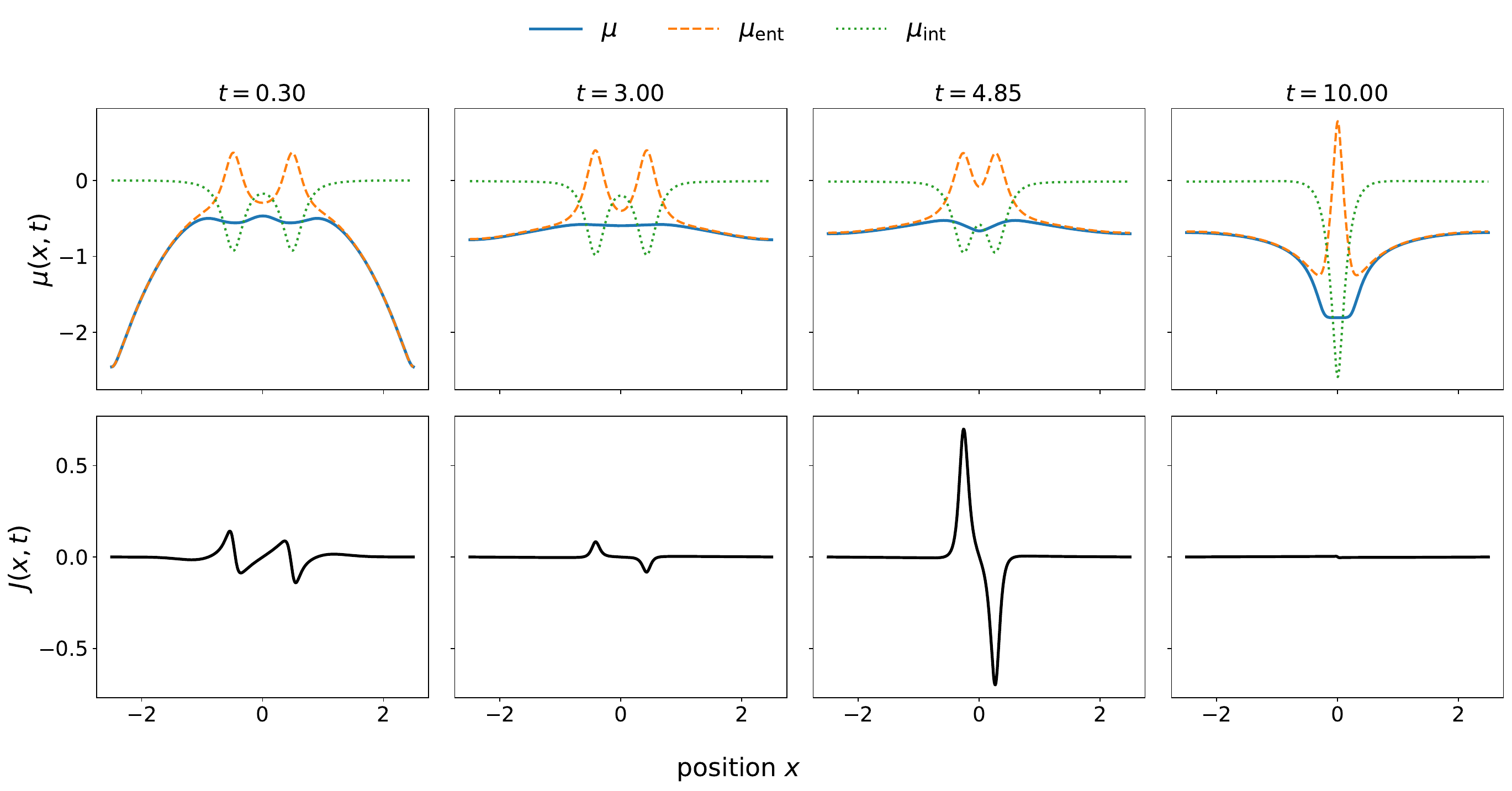}
\caption{\new{\textbf{Local transport mechanisms during the non-monotone transient of~\Cref{ex:twoGaussians}.}
Chemical potential $\mu=\frac{\delta \mathcal{F}}{\delta \rho}$ (top row) and associated flux $J=-\rho\nabla\mu$ (bottom row) at four representative times. The chemical potential is decomposed as $\mu=\mu_{\rm ent}+\mu_{\rm int}$. At $t=0.3$, the flux is directed towards the two density peaks, corresponding to their initial growth; at $t=3$, during the diffusion-dominated stage, the flux is substantially weaker while the two clusters approach each other; at $t=4.85$, a pronounced inward flux towards the centre accompanies the merging of the two clusters; and at $t=10$, close to the final clustered state, the remaining flux is small.}
} 
\label{fig:twoGaussians_2}
\end{figure}

\begin{ex}[Alternating transient regimes relaxing to a clustered state] 
\label{ex:twoGaussians}
We consider the same Morse-type interaction potential and parameter values as in~\Cref{ex1}, but reduce the noise strength to
\[
\sigma=0.65,
\]
which still lies in the coexistence regime. As initial condition we choose the normalized superposition of two periodized Gaussian profiles,
\[
\rho_0
\sim
\tfrac12\Big(
\mathcal N_{\mathrm{per}}(0.5,0.2^2)
+
\mathcal N_{\mathrm{per}}(-0.5,0.2^2)
\Big).
\]
    
\begin{figure}
 \begin{subfigure}{0.41\textwidth}
        \centering
    \includegraphics[width=\textwidth]{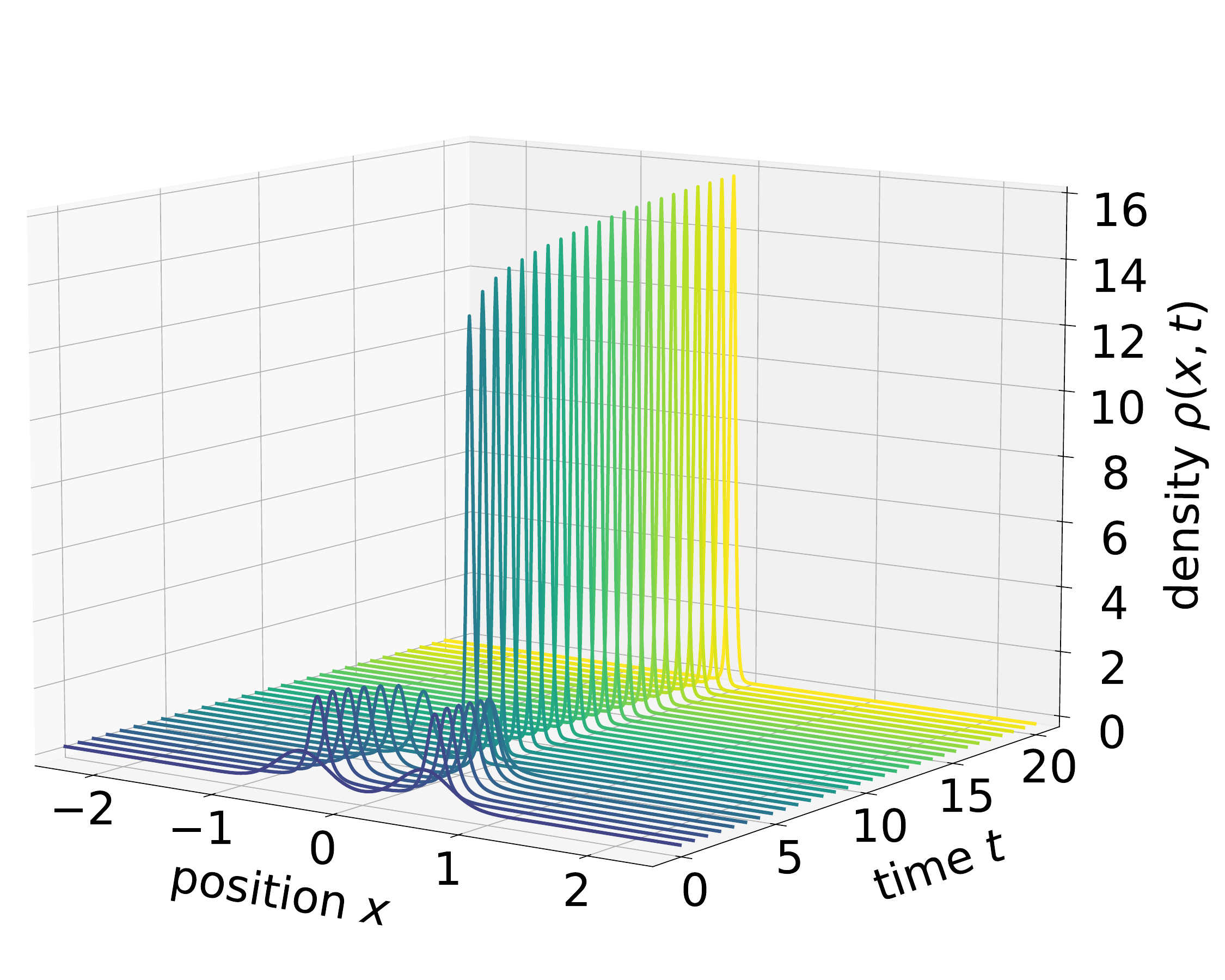}
     \caption{Density evolution $\rho(x,t)$}
    \end{subfigure}
     \begin{subfigure}{0.58\textwidth}
        \centering
    \includegraphics[width=\textwidth]{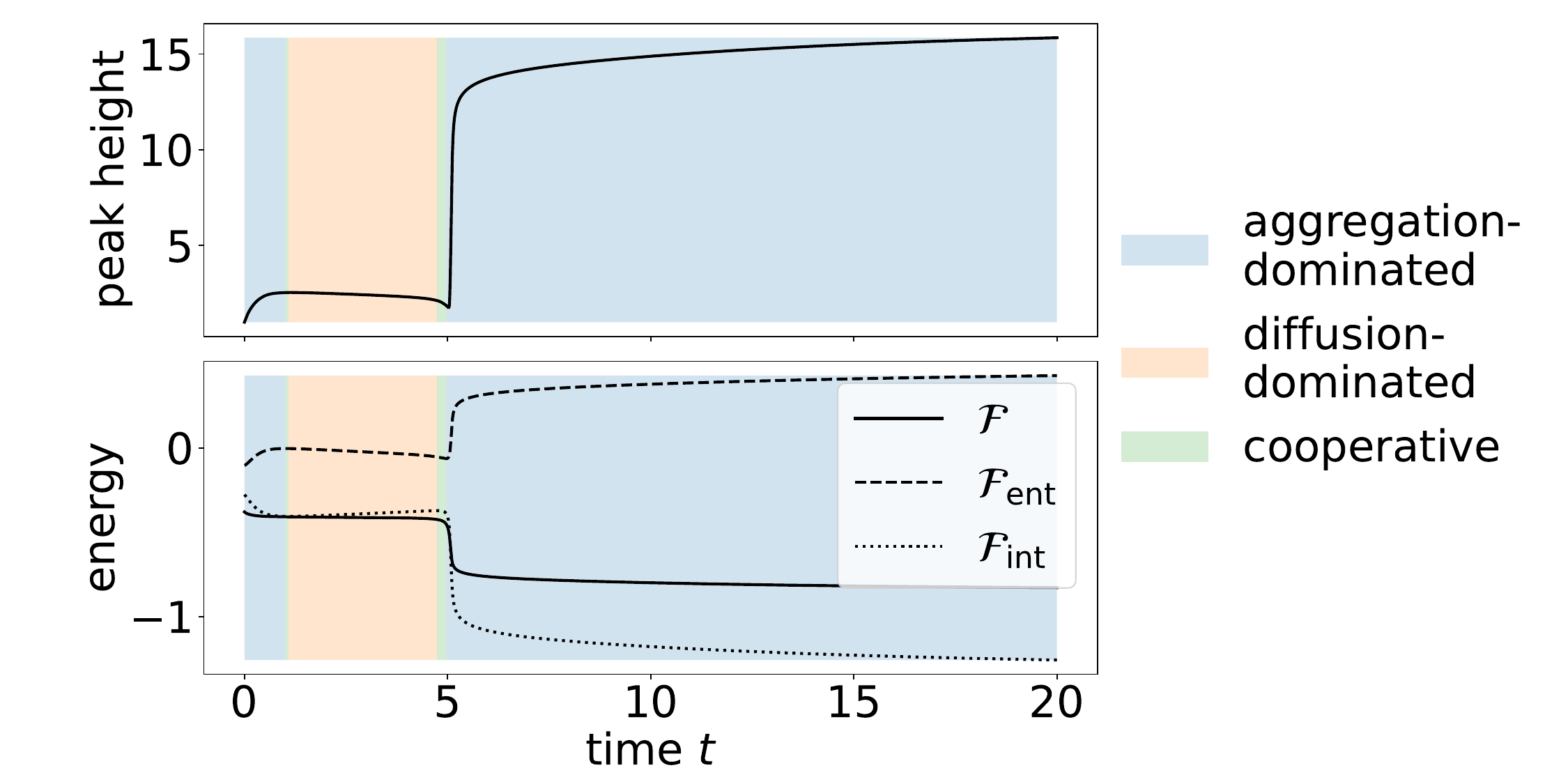}
     \caption{Peak height and free-energy contributions}
    \end{subfigure}
    \caption{\textbf{Alternating transient regimes during relaxation to the clustered equilibrium.} 
    Numerical solution of the PDE~\eqref{eq:MKV} for $\sigma= 0.65$ starting from~$\rho_0 \sim \tfrac{1}{2} (\mathcal{N}_{\text{per}}(-0.5,0.2^2)+\mathcal{N}_{\text{per}}(0.5,0.2^2))$, see~\Cref{ex:twoGaussians}. 
    (a) Time evolution of the spatial density profiles  $\rho(x,t)$, showing the formation, approach, and eventual merging of two clusters.
    (b) Time evolution of the peak height and free-energy contributions. The peak height exhibits alternating increasing and decreasing phases, while the free-energy contributions reveal successive aggregation-dominated, diffusion-dominated, cooperative, and renewed aggregation-dominated regimes.
    }
    \label{fig:twoGaussians}
\end{figure} 

The resulting dynamics \new{are illustrated in~\Cref{fig:twoGaussians}. As in~\Cref{ex1},} they exhibit alternating transient behaviour, now converging towards a clustered stationary state. Initially, both peaks increase in height while the interaction energy decreases and the entropy increases, indicating an aggregation-dominated regime.

At intermediate times ($1.0\lesssim t\lesssim4.7$), the two clusters slowly approach each other while their peak heights decrease. During this phase, the interaction energy increases and the entropic contribution decreases, corresponding to a diffusion-dominated regime. Despite the decreasing peak heights, the clusters continue to move towards each other due to the spatial structure of the interaction field. \new{This is followed by a short cooperative relaxation regime ($4.7\lesssim t\lesssim5.0$), during which both energy contributions decrease simultaneously.} \sout{A short interval of cooperative relaxation follows ($4.7\lesssim t\lesssim5.0$), during which both energy contributions decrease simultaneously.}

At later times ($t\gtrsim5.0$), the two clusters merge into a single peak. After merging, the interaction energy decreases again while the entropy increases, marking a renewed aggregation-dominated regime. The resulting cluster sharpens and converges towards a stable stationary clustered state.

\new{In analogy to~\Cref{fig:oneGaussian_2}, the corresponding chemical-potential and flux profiles are shown in~\Cref{fig:twoGaussians_2}, illustrating how the local transport structure changes during the approach and merging of the two clusters and how the flux becomes small near the stationary state.}

This example demonstrates that alternating transient regimes may arise even when the asymptotic state is clustered, showing that transient diffusion-dominated dynamics can play an essential role in the formation and merging of clusters. \hfill$\diamond$
\end{ex}

\new{\paragraph{Geometric observables versus regime classification.}
Taken together, the two examples show that geometric clustering observables
and the free-energy-based classification
 describe different aspects of the
transient dynamics. 
Changes in peak height may accompany regime changes, but the correspondence is not one-to-one. 
In particular, local concentration may increase even while the evolution is globally diffusion-dominated, while changes in cluster position or merging may occur without being captured by peak height alone. 
The sensitivity of these behaviours with respect to the initial condition is further illustrated in Appendix~\ref{app:initial_state}.}

\section{Conclusion} 
\label{sec:conclusion}

In this work, we investigated transient clustering dynamics in a McKean--Vlasov aggregation--diffusion model with locally attractive interactions on the periodic domain. Using numerical solutions of the PDE for varying noise strengths, we focused on the transient evolution preceding convergence towards different types of stationary states. 
Starting from non-uniform initial conditions, we observed alternating aggregation- and diffusion-dominated regimes, together with non-monotone evolution of geometric observables such as the density peak height. \new{Using the temporal evolution of the entropic and interaction contributions to the free energy as a diagnostic, we found that the dominant mechanism can switch repeatedly along a single relaxation trajectory, despite monotone dissipation of the total free energy.}
In particular, we studied two representative scenarios within the coexistence regime: relaxation towards the homogeneous equilibrium and convergence towards a clustered stationary state. \sout{In both cases, the transient dynamics exhibited alternating energetic regimes despite monotone dissipation of the total free energy.}

Our analysis highlights that different observables capture different aspects of clustering dynamics and may therefore exhibit different forms of non-monotone behaviour. In particular, local geometric quantities such as peak height are related to, but not uniquely determined by, the underlying energetic mechanisms governing the evolution. From this perspective, monotonicity should not necessarily be viewed as a defining criterion for a useful \new{observable that monitors clustering dynamics}. Rather, non-monotone behaviour may itself encode relevant transient transport mechanisms.
The energetic perspective based on the free-energy decomposition provides a natural variational and thermodynamic framework for distinguishing aggregation- and diffusion-dominated transient regimes. 

In our numerical examples, the observed alternating behaviour arises for noise strengths close to critical values and depends sensitively on the initial condition, reflecting the complex geometry of the underlying free-energy landscape rather than solely the equilibrium structure of the system. \new{Overall, the results demonstrate that geometric clustering observables and the analysis of the free-energy contributions provide complementary information about transient clustering dynamics and the underlying transport mechanisms in interacting particle systems.}

Several directions for future work remain open. Beyond the representative examples considered here, it would be interesting to investigate more systematically for which interaction kernels, noise strengths, and initial conditions alternating transient regimes emerge.
Another natural direction concerns the relation between the free-energy decomposition and alternative observables for quantifying clustering dynamics. Possible candidates include moment-based quantities, Fourier modes, or data-driven collective variables. 
From an application-oriented perspective, distinguishing aggregation- and diffusion-dominated regimes may also provide insight into transient pattern formation in systems arising from swarming dynamics or opinion formation models.
Finally, it remains an open question whether similar forms of transient non-monotone behaviour may occur outside the coexistence regime, in particular in parameter regimes with a unique stable stationary state.

\paragraph{Acknowledgment} This research was supported by the Deutsche Forschungsgemeinschaft (DFG) through the Collaborative Research Center CRC 1114 ``Scaling Cascades in Complex Systems'' (Project No.~235221301), in particular through collaboration between project C03 ``Multiscale modelling and simulation for spatiotemporal master equations'' and project A05 ``Probing Scales in Equilibrated Systems by Optimal Nonequilibrium Forcing''. Additional support was provided through Germany's Excellence Strategy -- MATH+: Berlin Mathematics Research Center (EXC 2046/1, project 390685689). C.H.\ was further supported by the German Federal Government, the Federal Ministry of Research, Technology and Space, and the State of Brandenburg within the framework of the joint project EIZ: Energy Innovation Center (project numbers 85056897 and 03SF0693A) with funds from the Structural Development Act (Strukturstärkungsgesetz) for coal-mining regions.

\paragraph{Use of generative AI} During the revision of this manuscript, the authors used ChatGPT (OpenAI, GPT-5.6, October 2026) to assist with language editing, rephrasing, and drafting parts of the text. 
All AI-generated suggestions were critically reviewed and edited by the authors, who take full responsibility for the content of the manuscript.

\paragraph{Competing interests.} The authors declare no competing interests.

\appendix

\section{Appendix}

\begin{figure}[t]
\centering
\includegraphics[width=0.5\textwidth]{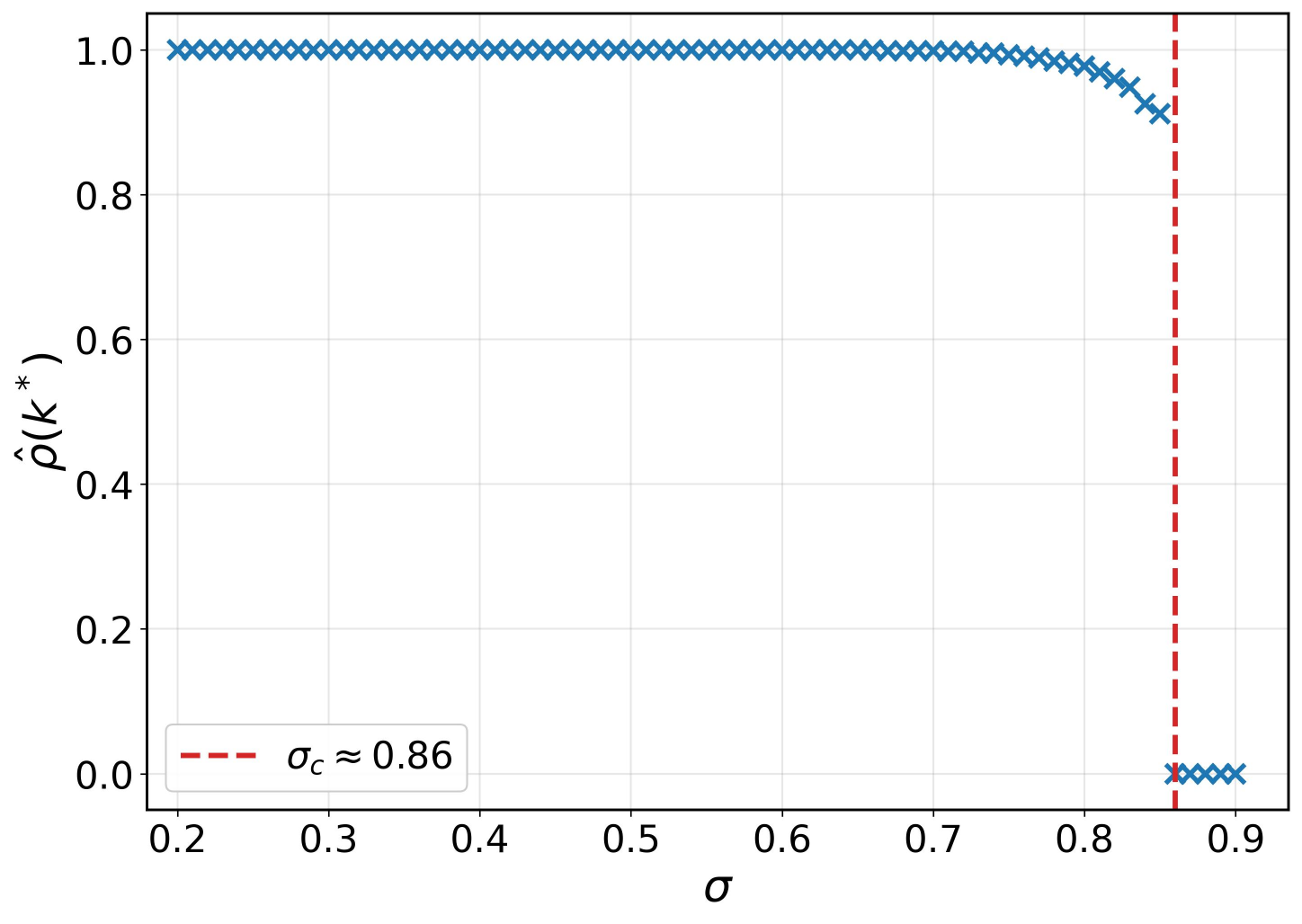}
\caption{\textbf{Numerical estimation of the critical noise strength $\sigma_c$.}
The modulation amplitude $A(\rho_T)=|\hat{\rho}_T(k^*)|$ of the density at time $T=100$
is shown for different values of $\sigma$, starting from the strongly
clustered initial condition
$\rho_0 \sim \mathcal{N}_{\mathrm{per}}(0,0.1^2)$.
The abrupt transition from positive values to values close to zero
occurs at approximately $\sigma_c \approx 0.86$.}
\label{fig:sigmac_estimation}
\end{figure}

\new{\subsection{Numerical estimation of the critical noise strength $\sigma_c$}
\label{app:sigmac}
To estimate the critical noise strength $\sigma_c$ numerically, we solve the McKean--Vlasov equation for increasing values of $\sigma$, starting from strongly clustered initial densities, and monitor the modulation amplitude of the final density,
\begin{equation}
    A(\rho_T)
    :=
    |\hat{\rho}_T(k^*)|
    =
    \left|
    \frac{1}{L}
    \int_{-L/2}^{L/2}
    \rho_T(x)
    e^{-2\pi i k^*x/L}\,dx
    \right|,
\end{equation}
where $k^* \in \mathbb{N}_+$ denotes the dominant Fourier mode associated with the clustered state. 
The simulations are run up to $T=100$. For the tested clustered initial conditions, $A(\rho_T)$ remains clearly positive for smaller noise strengths and drops abruptly to values close to zero when $\sigma$ exceeds approximately $0.86$. We therefore use $\sigma_c \approx 0.86$ as a numerical estimate for the upper boundary of the coexistence regime. Figure~\ref{fig:sigmac_estimation} 
illustrates this abrupt transition for the initial condition $ \rho_0 \sim \mathcal{N}_{\mathrm{per}}(0,0.1^2)$.}

\subsection{Sensitivity with respect to the initial state}
\label{app:initial_state} 

\begin{figure}
 \begin{subfigure}{0.41\textwidth}
        \centering
    \includegraphics[width=\textwidth]{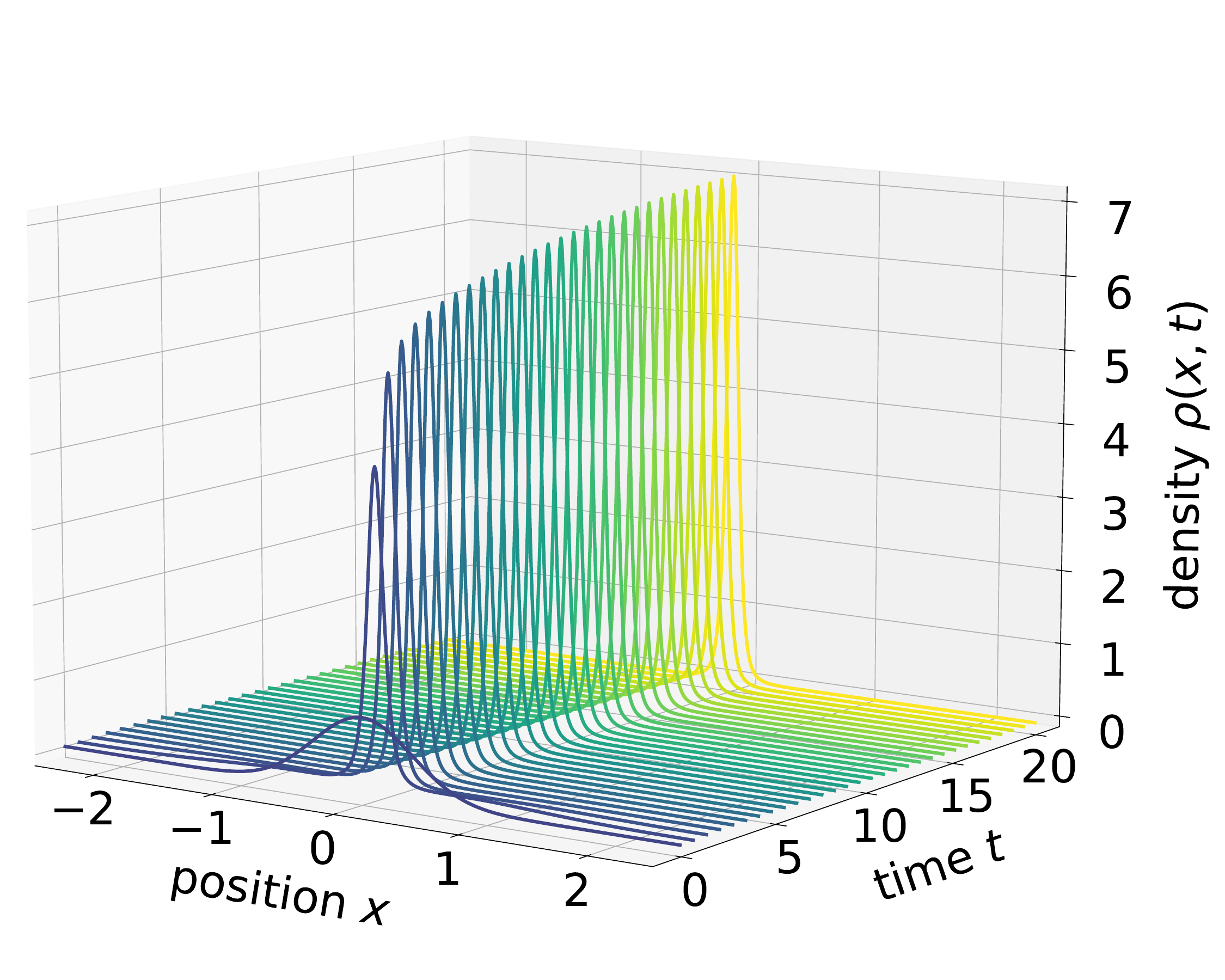}
     \caption{Density evolution $\rho(x,t)$}
    \end{subfigure}
         \begin{subfigure}{0.58\textwidth}
        \centering
    \includegraphics[width=\textwidth]{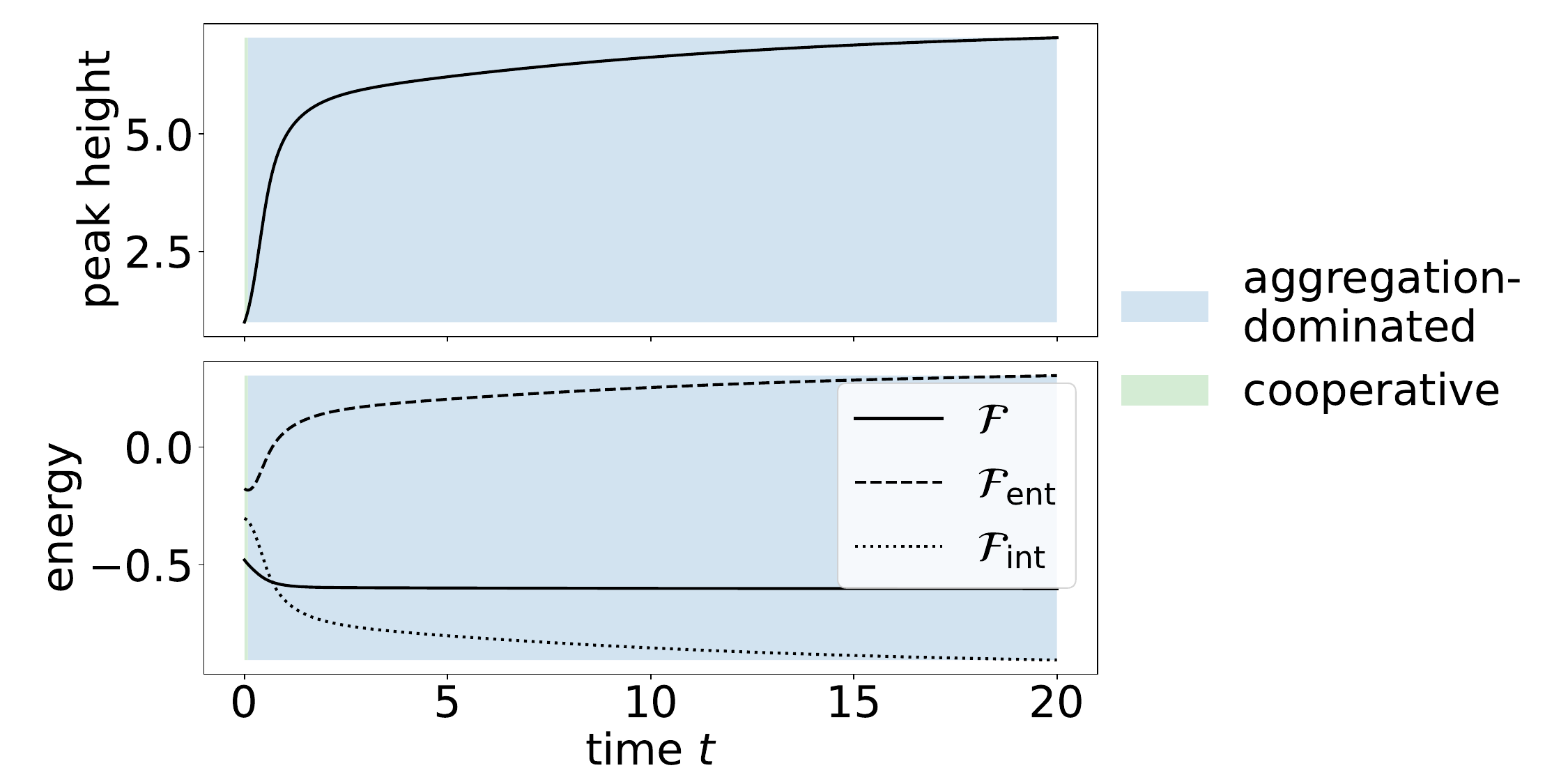}
     \caption{Peak height and free-energy contributions}
    \end{subfigure}
    \caption{\textbf{Aggregation-dominated relaxation to the clustered equilibrium.} Numerical solution of the PDE~\eqref{eq:MKV} for $\sigma= 0.838$ starting from~$\rho_0\sim \mathcal{N}_{\text{per}}(0,0.4^2)$. 
   (a) Time evolution of the spatial density profiles~$\rho(x,t)$, showing concentration towards a clustered equilibrium.
    (b) Time evolution of the peak height and free-energy contributions. The peak height increases monotonically, while the free-energy contributions indicate predominantly aggregation-dominated dynamics, preceded by a short cooperative regime.
   }
    \label{fig:oneGaussian_start0.4}
\end{figure}

Both the transient dynamics and the long-time limit depend sensitively on the initial condition. Keeping all parameters identical to those in~\Cref{ex1} from \Cref{sec:alternating} while varying only the initial state~$\rho_0$, we observe that choosing $\rho_0 \sim \mathcal{N}_{\text{per}}(0,0.4^2)$, 
corresponding to a slightly smaller initial variance, leads to predominantly aggregation-dominated dynamics converging to a single-peak clustered equilibrium; see~\Cref{fig:oneGaussian_start0.4}. 

\begin{figure}
 \begin{subfigure}{0.41\textwidth}
        \centering
    \includegraphics[width=\textwidth]{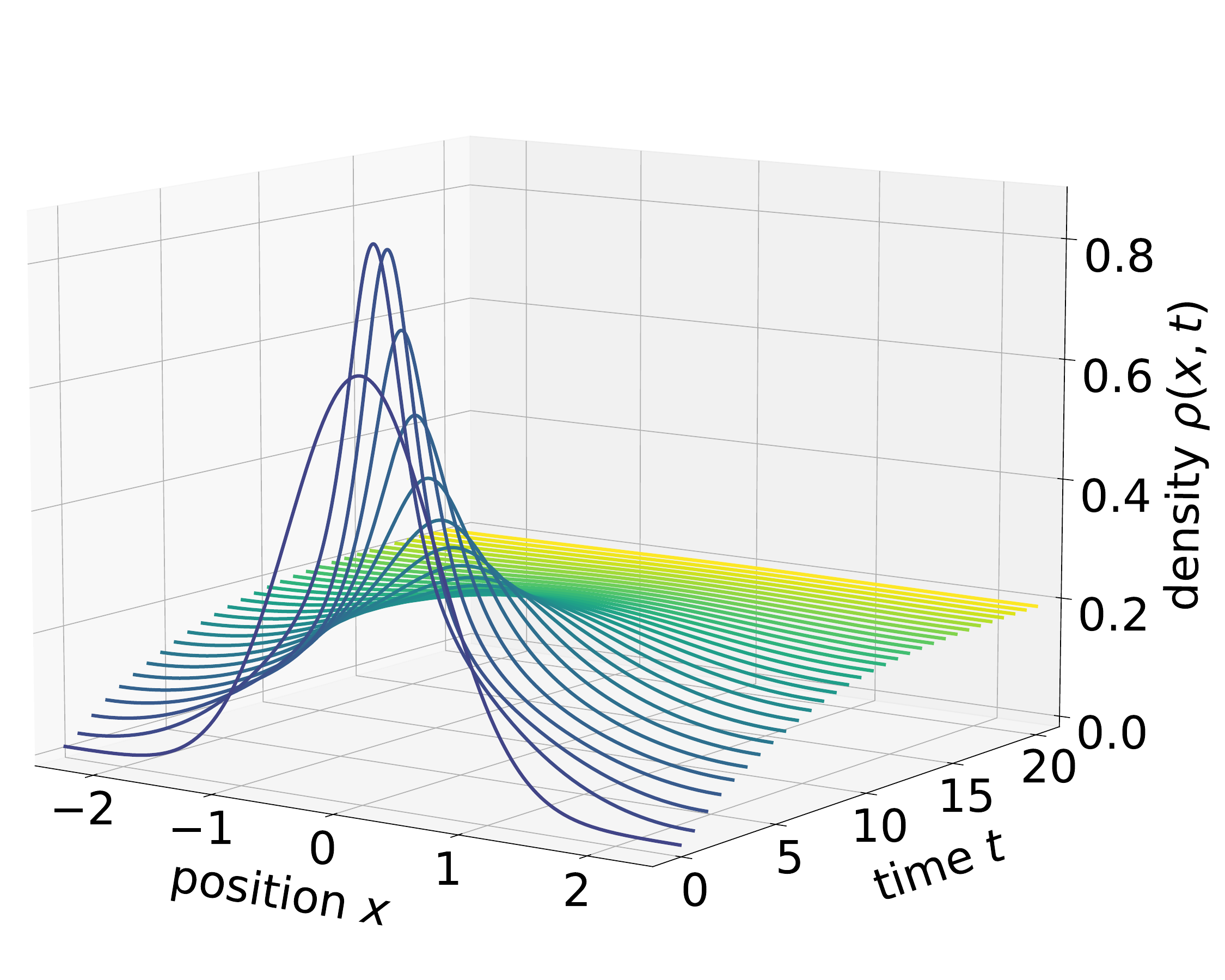}
     \caption{Density evolution $\rho(x,t)$}
    \end{subfigure}
         \begin{subfigure}{0.58\textwidth}
        \centering
    \includegraphics[width=\textwidth]{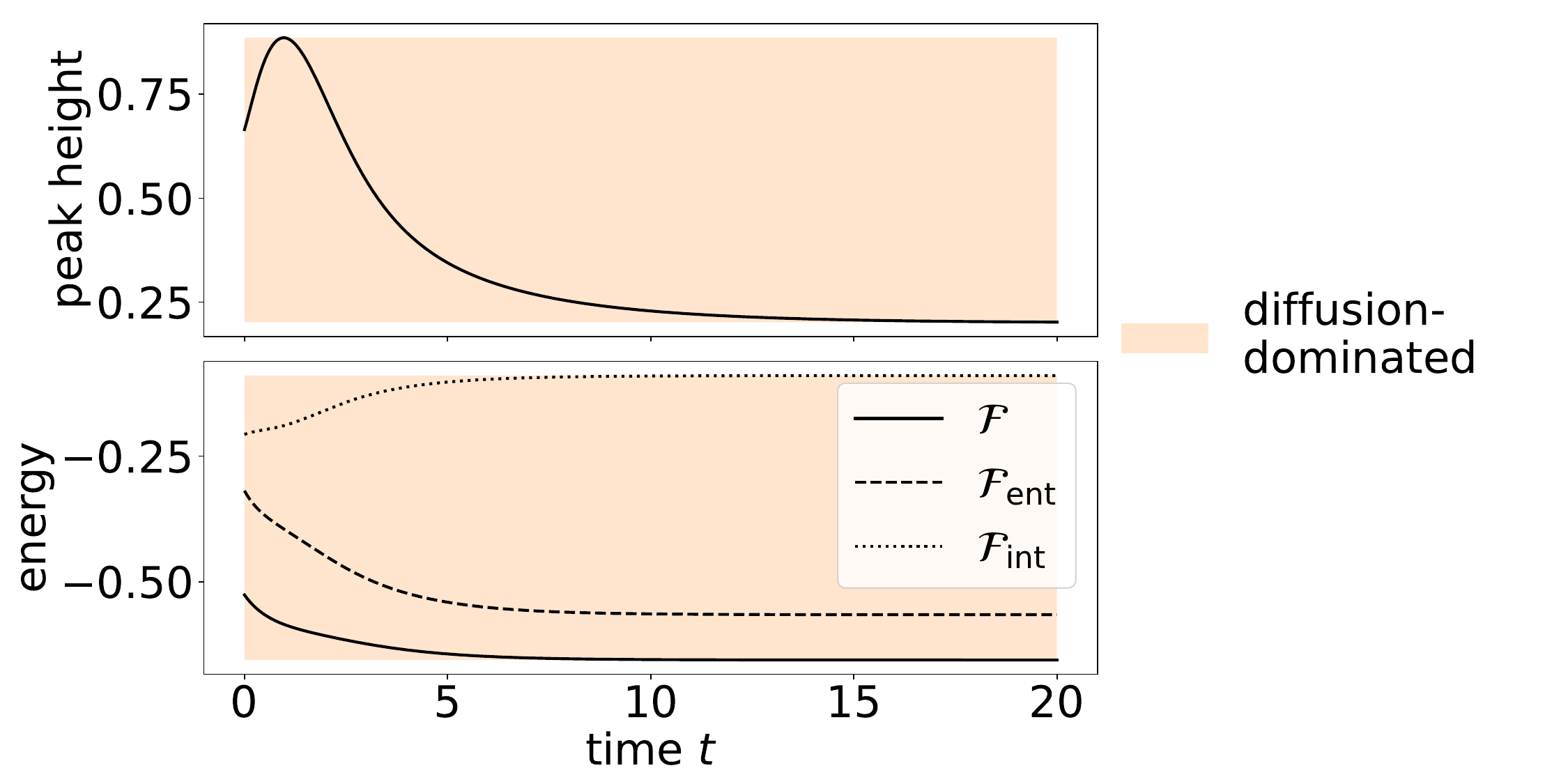}
     \caption{Peak height and free-energy contributions}
    \end{subfigure}
    \caption{\textbf{Diffusion-dominated relaxation to the homogeneous equilibrium (but non-monotonicity in peak height).} Numerical solution of the PDE~\eqref{eq:MKV} for $\sigma= 0.838$ starting from~$\rho_0\sim \mathcal{N}_{\text{per}}(0,0.6^2)$. 
    (a) Time evolution of the spatial density profiles~$\rho(x,t)$, showing initial concentration followed by spreading towards the homogeneous equilibrium.
    (b) Time evolution of the peak height and free-energy contributions. The peak height first increases and subsequently decreases, while the free-energy contributions indicate diffusion-dominated dynamics throughout the evolution.    
    }
    \label{fig:oneGaussian_start0.6}
\end{figure}

In contrast, increasing the initial variance such that $\rho_0 \sim \mathcal{N}_{\text{per}}(0,0.6^2)$ results in purely diffusion-dominated dynamics converging to the homogeneous equilibrium; see~\Cref{fig:oneGaussian_start0.6}.
Nevertheless, the peak height still exhibits a temporary increase before decreasing again, showing that the maximal density alone is not a reliable indicator of aggregation- or diffusion-dominated behaviour. \new{This further motivates complementing geometric clustering observables by an energetic analysis of the transient dynamics.}

\subsection{Choice and number of observables}
\label{sec:choice_and_number_of_observables}

The numerical experiments reveal that peak height is, in general, not a good ``progress variable'' for unambiguously describing the relaxation to equilibrium, since it may evolve non-monotonically in time. \new{Nevertheless, its intuitive association with the concentration of mass makes it a natural geometric observable for monitoring clustering. The results of Section \ref{subsec:Numerical investigation} illustrate, however, the limitations of interpreting peak height as a measure of clustering strength.}

\new{Another observable that quantifies the spread of a density is its second moment
\begin{equation} 
\label{eq:second_moment}
     \int x^2 \rho(x,t) \, dx.
\end{equation}
For an initial density $\rho_0$ that is symmetric around the origin and an even interaction potential $U$, the symmetry of the PDE solution is preserved for all
times and its mean remains zero. In this case, the second moment coincides with the variance.}

\begin{figure}
\begin{subfigure}{0.31\textwidth}
    \centering
\includegraphics[width=\textwidth]{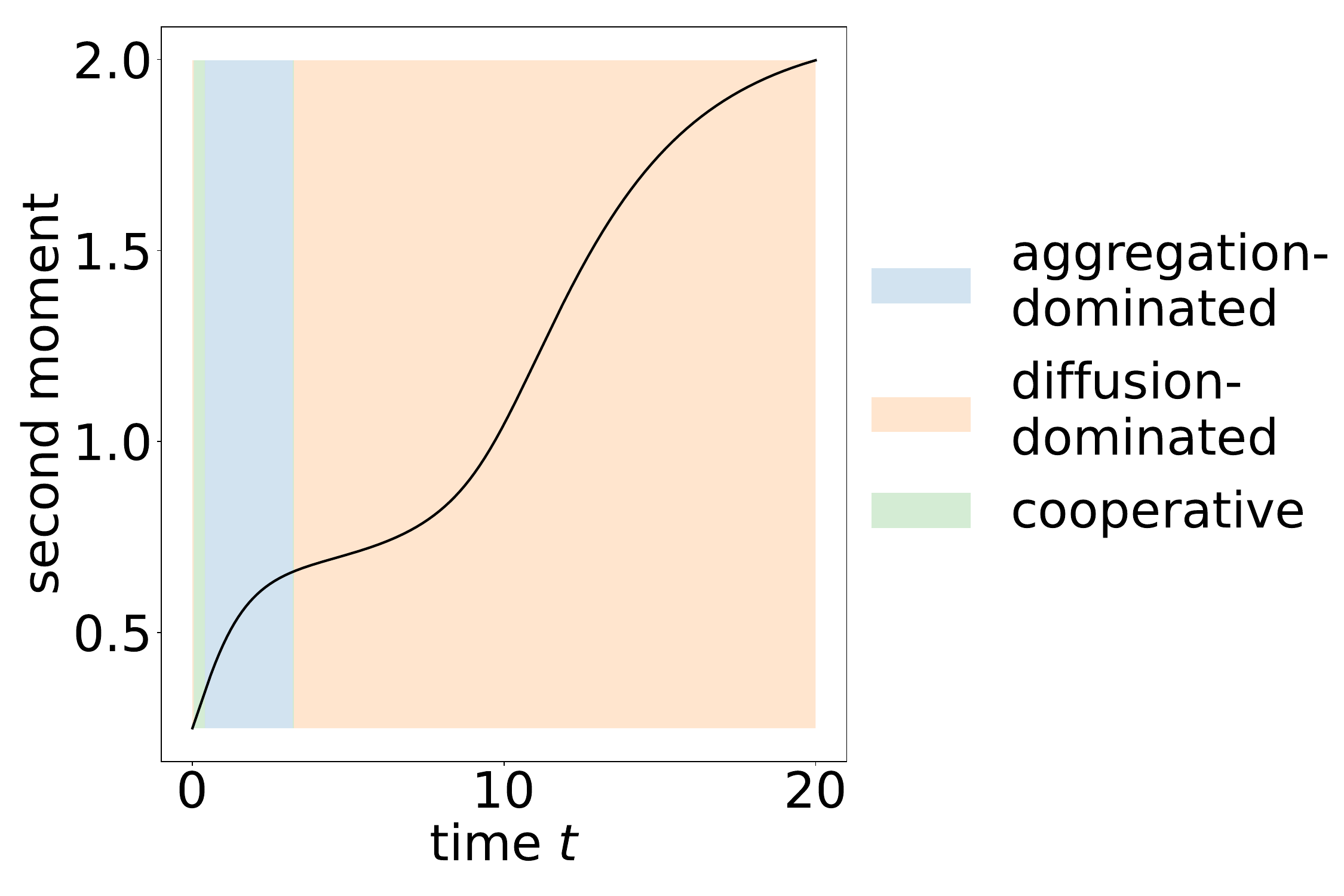}
 \caption{Second moment \newline}
 \label{fig:variance_a}
\end{subfigure}
        \begin{subfigure}{0.32\textwidth}
    \centering
\includegraphics[width=\textwidth]{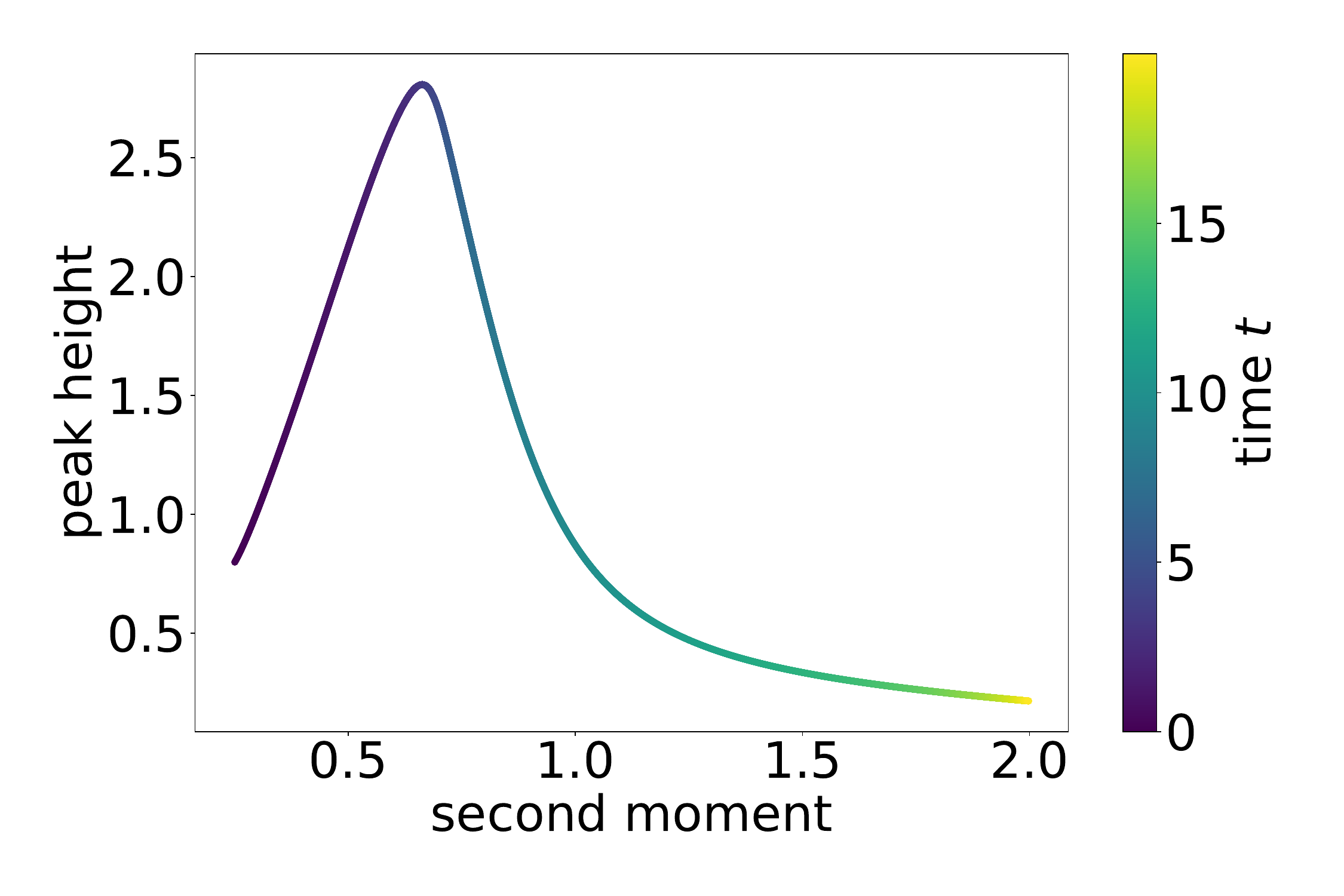}
 \caption{Second moment vs.\ peak height \newline coloured by time}
  \label{fig:variance_b}
\end{subfigure}
        \begin{subfigure}{0.32\textwidth}
    \centering
\includegraphics[width=\textwidth]{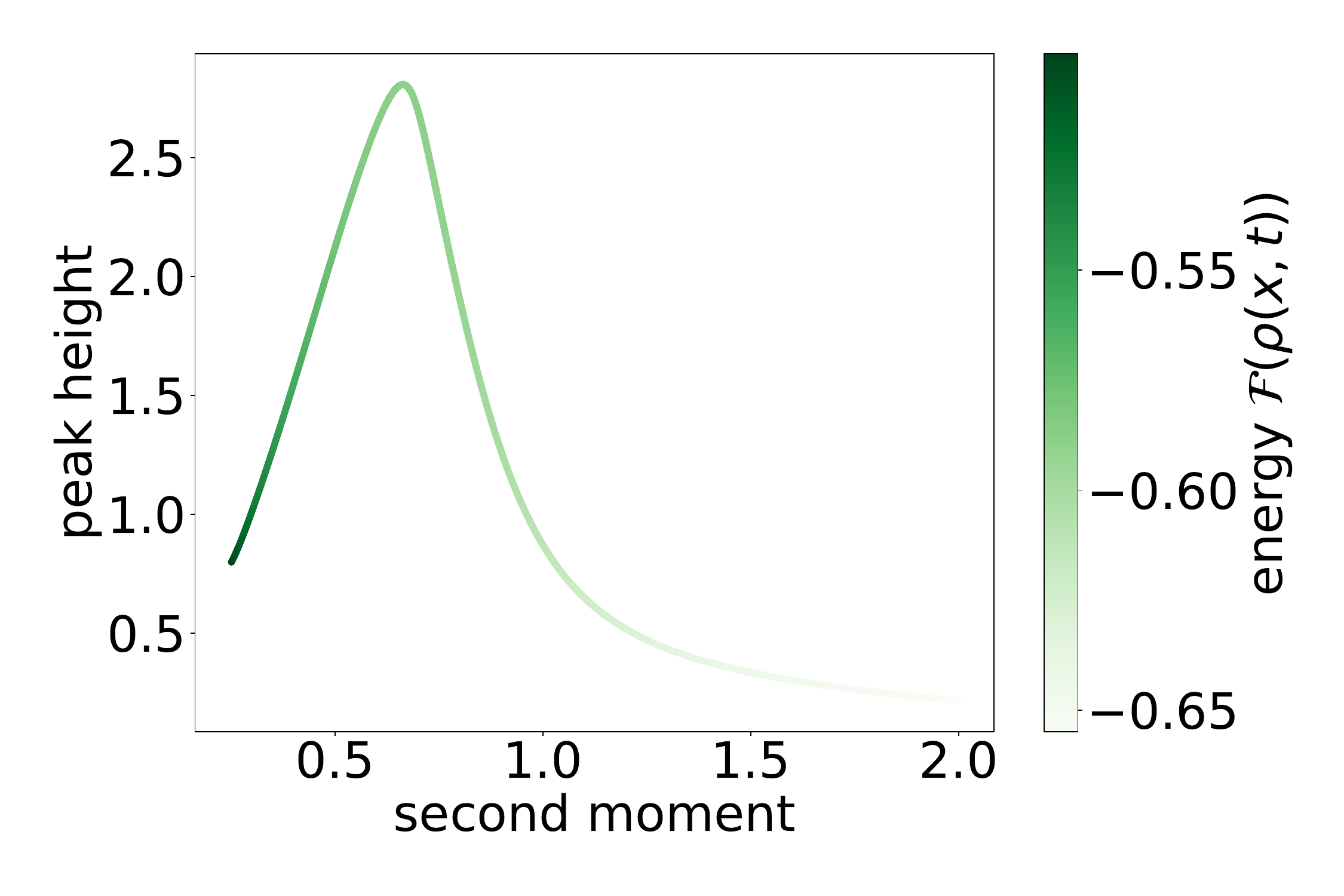}
 \caption{Second moment vs.\ peak height \newline coloured by free energy}
  \label{fig:variance_c}
\end{subfigure}
\caption{\textbf{Second moment, peak height, and free energy for~\Cref{ex1}.}
(a) Temporal evolution of the second moment, with background colours indicating
the energetic regimes.
(b) Trajectory in the plane spanned by the second moment and peak height, with colour indicating time.
(c) The same trajectory, with colour indicating the total free energy
$\mathcal{F}[\rho(\cdot,t)]$.}
\label{fig:variance}
\end{figure}

For the setting of~\Cref{ex1}, the second moment turns out to be a useful progress variable in the sense that it is monotone in time; see~\Cref{fig:variance_a}. At the same time, however, it does not clearly separate aggregation-, diffusion-, and cooperative regimes, since the curve exhibits similar features across the different phases. This suggests that a single observable may not suffice to identify the energetic regimes and that higher-dimensional descriptions may be required.

Indeed, the middle and right panels of~\Cref{fig:variance} show the relaxation dynamics as a function of both second moment \textit{and} peak height. Comparing the left and middle panels, we observe that the local maximum in peak height roughly corresponds to the regime transition. In addition, the right panel demonstrates that the free energy changes only weakly during the regime switch, providing further evidence for the interplay between entropy- and interaction-driven exploration of the free-energy landscape.

\subsection{Hegselmann--Krause model}
\label{sec:Hegselmann--Krause model}

\begin{figure}
     \begin{subfigure}{0.41\textwidth}
        \centering
    \includegraphics[width=\textwidth]{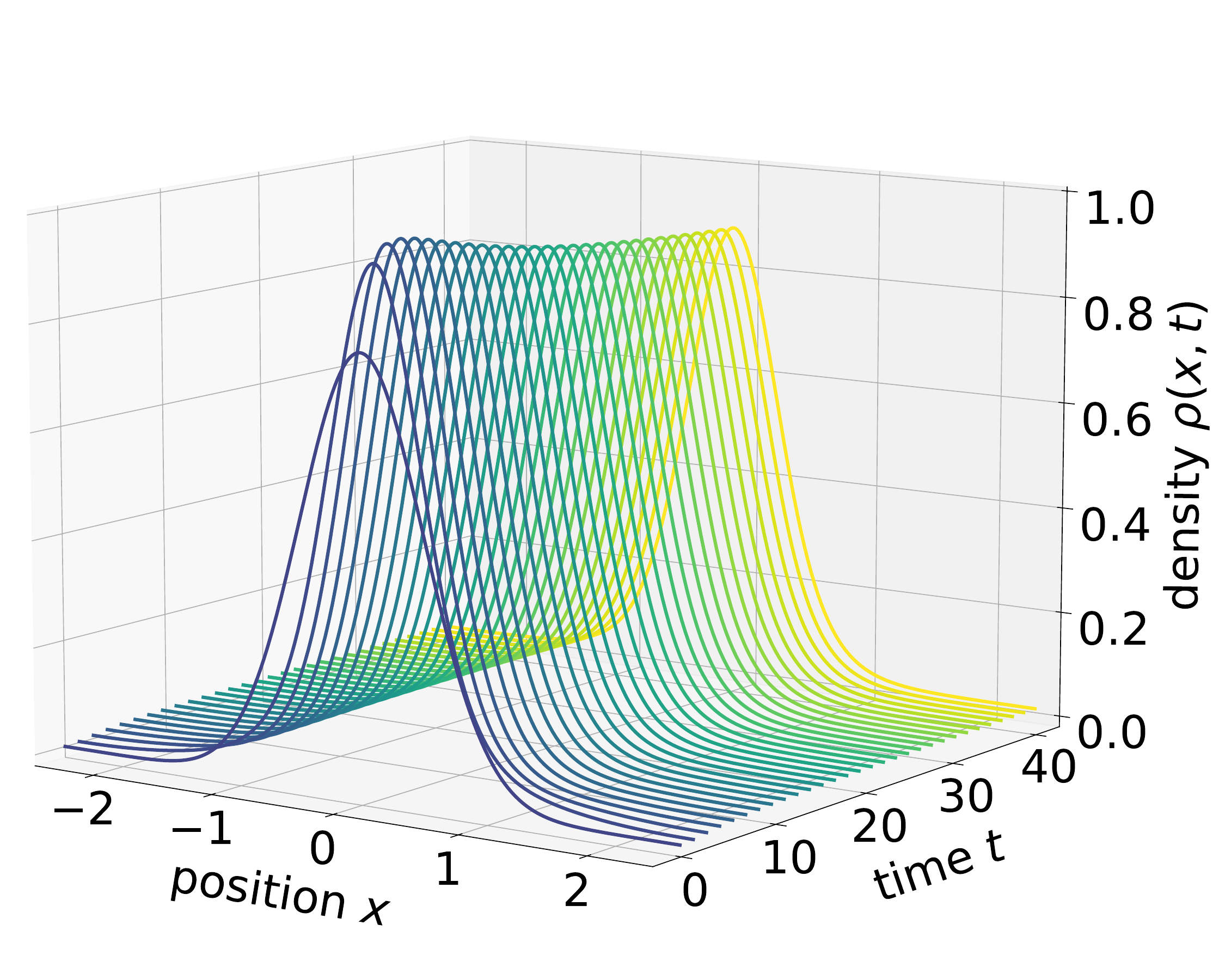}
     \caption{Density evolution $\rho(x,t)$}
    \end{subfigure}
    \begin{subfigure}{0.58\textwidth}
        \centering
    \includegraphics[width=\textwidth]{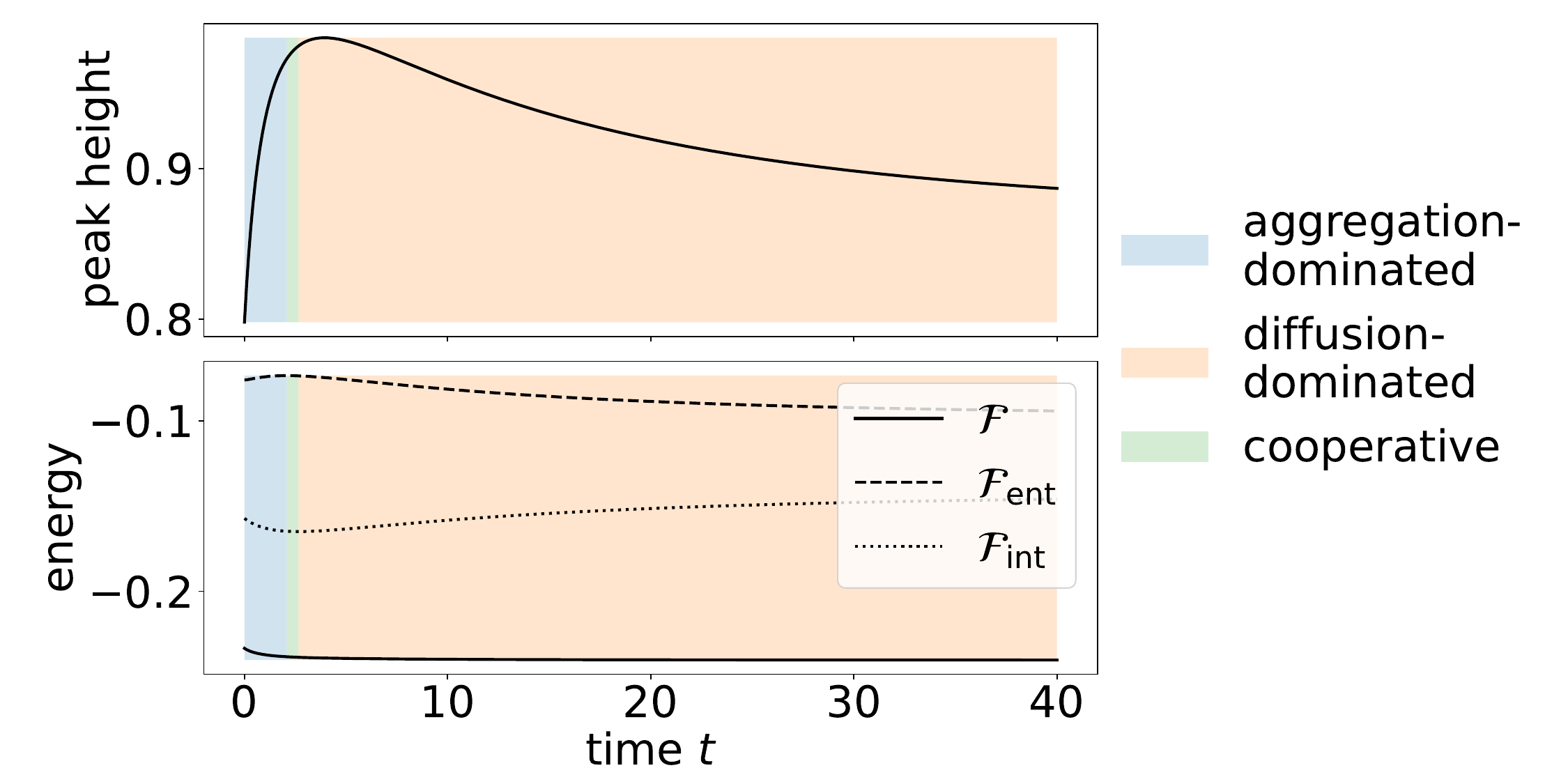}
     \caption{Peak height and free-energy contributions}
    \end{subfigure}
    \caption{\textbf{Hegselmann--Krause: Alternating transient regimes during relaxation towards the homogeneous equilibrium.} Numerical solution of the PDE~\eqref{eq:MKV} for $U(x)=\frac{1}{2}(|x|^2-R_0^2)\mathbbm{1}_{[0,R_0]}(|x|)$ with $R_0=1$ and $\sigma= 0.458$ starting from~$\rho_0\sim \mathcal{N}_{\text{per}}(0,0.5^2)$. 
    (a) Time evolution of the spatial density profiles~$\rho(x,t)$, showing initial concentration followed by spreading towards the homogeneous equilibrium.
    (b) Time evolution of the peak height and free-energy contributions. The peak height first increases and subsequently decreases, while the free-energy contributions reveal successive aggregation-dominated, cooperative, and diffusion-dominated regimes.
    }
    \label{fig:oneGaussian_HK}
    \end{figure}

Similar to \Cref{ex1} in \Cref{sec:alternating}, we investigate another example of a locally attractive interaction potential---distinct from the Morse potential with the parameter values given in \Cref{eq:parameter_values}---that exhibits alternating transient regimes. The \textit{Hegselmann--Krause model}, given in the mean-field limit by \Cref{eq:MKV} with interaction kernel $U(x)=\frac{1}{2}(|x|^2-R_0^2)\mathbbm{1}_{[0,R_0]}(|x|)$, where $R_0>0$, describes aggregation phenomena in opinion dynamics \citep{hegselmann2002, garnier2017consensus, gerber2025formation}. Analogous to the Morse potential, the truncated quadratic potential generates equilibria distinct from the homogeneous state for sufficiently small noise strengths~\citep{carrillo_long-time_2020}.
Figure~\ref{fig:oneGaussian_HK} shows that, for a suitable value of the noise strength $\sigma$ (within the phase coexistence regime), alternating transient regimes arise during relaxation to the homogeneous equilibrium. Over the time interval shown, the density has not yet reached the homogeneous equilibrium; longer simulations confirm continued relaxation towards the uniform state.

\bibliographystyle{agsm}
\bibliography{references}

\end{document}